\documentclass[journal,11pt,draftcls,onecolumn]{IEEEtran}

\newcommand{\revised}[1]{{\color{black} #1}}

\usepackage{color}
\usepackage{cite}

\usepackage{clrscode}
\usepackage{amsmath, amssymb, amsfonts, amsthm}
\usepackage{epsfig, epstopdf}
\usepackage{enumerate}
\usepackage{dsfont}
\usepackage{graphicx}
\usepackage{subfig}

\theoremstyle{theorem}
\newtheorem{thm}{Theorem}[]
\newtheorem{coro}{Corollary}[]
\newtheorem{lem}{Lemma}

\theoremstyle{remark}
\newtheorem{rem}{Remark}
\theoremstyle{definition}
\newtheorem{defn}{Definition}

\newcommand{\subg}[1]{{\partial}{#1}}
\newcommand{\prox}[1]{{P}_{#1}}

\newcommand{\Lp}[3]{{\left\|{#1}\right\|}_{#2}^{#3}}
\newcommand{\LCr}[1]{\left\{#1\right\}}
\newcommand{\LPr}[1]{\left(#1\right)}

\newcommand{\RR}[1]{\mathbb{R}^{#1}}

\newcommand{\amin}[1]{\underset{#1}{\operatorname{argmin}}}

\newcommand{\mat}[1]{{\bf{#1}}}
\newcommand{\mats}[1]{{\boldsymbol{#1}}}

\newcommand{\identi}{\overset{\small \Delta}{\operatorname{=}}}

\renewcommand{\identi}{\triangleq}
\newcommand{\vstack}{\mathop{\mathrm{vec}}\nolimits}

\title{\huge Nested Sparse Approximation: Structured Estimation of V2V Channels Using
Geometry-Based Stochastic Channel Model}
\author{Sajjad~Beygi~\IEEEmembership{Student Member, IEEE},
 Urbashi~Mitra~\IEEEmembership{Fellow, IEEE}, and Erik~G.~Str\"om$^{*}$~\IEEEmembership{Senior Member, IEEE}
 
 \thanks{This research was funded in part by one or all of these grants:  
ONR N00014-09-1-0700,
AFOSR FA9550-12-1-0215, 
DOT CA-26-7084-00, 
NSF CCF-1117896,  
NSF CNS-1213128,  
NSF CCF-1410009,  
NSF CPS-1446901, 
Barbro Osher Pro Suecia Foundation, 
Ericsson's Research Foundation FOSTIFT-13:038, 
and Adlerbert Research Foundation.}
\thanks{School of Electrical Engineering,
University of Southern California
$^{*}$Dept. of Signals and Systems,
Chalmers University of Technology
Emails: \{beygihar, ubli\}@usc.edu,  erik.strom@chalmers.se}}

\begin{document}
\maketitle
\IEEEpeerreviewmaketitle

\begin{abstract}
Future intelligent transportation systems promise increased safety and energy efficiency.  Realization of such systems will require vehicle-to-vehicle (V2V) communication technology.  High fidelity V2V communication is, in turn, dependent on accurate V2V channel estimation.  V2V channels have characteristics differing from classical cellular communication channels.  Herein, geometry-based stochastic modeling is employed to develop a characterization of  V2V channel channels.  The resultant model exhibits significant structure; specifically, the V2V channel is characterized by three distinct regions within the delay-Doppler plane.  Each region has a unique combination of specular reflections and diffuse components resulting in a particular element-wise and group-wise sparsity.  This joint sparsity structure is exploited to develop a novel channel estimation algorithm.  A general machinery is provided to solve the jointly element/group sparse channel (signal) estimation problem using proximity operators of a broad class of regularizers.  The alternating direction method of multipliers using the proximity operator is adapted to optimize the mixed objective function.  
Key properties of the proposed objective functions are proven which ensure that the optimal solution is found by the new algorithm. 
The effects of pulse shape leakage are explicitly characterized and compensated, resulting in measurably improved performance.  Numerical simulation and real V2V  channel measurement data are used to evaluate the performance of the proposed method. Results show that the new method can achieve significant gains over previously proposed methods.  
\end{abstract}

\section{Introduction}
\IEEEPARstart{V}{ehicle}-to-vehicle (V2V) communication is central to future intelligent transportation systems, which will enable efficient and safer transportation with reduced fuel consumption \cite{matolak2013v2v}. In general,  V2V communication is anticipated to be short range with transmission ranges varying from a few meters to a few kilometers between two mobile vehicles on a road.  A big challenge of realizing V2V communication is the inherent fast channel variations (faster than in cellular \cite{bernado2013delay, nuckelt2011deterministic}) due to the mobility of both the transmitter and receiver.  Since channel state information can improve communication performance, a fast algorithm to accurately estimate V2V channels is of interest. Furthermore,  V2V channels are highly dependent on the geometry of road and the local physical environment \cite{molisch2009survey,matolak2013v2v}.  

A popular estimation strategy for fast time-varying channels is to apply Wiener filtering \cite{bello1963characterization, zemen2012adaptive}. Recently, \cite{zemen2012adaptive, bernado2013delay} present an adaptive Wiener filter to estimate  V2V channels using subspace selection. The main drawback of Wiener-filtering is that the knowledge of the scattering function is required \cite{bello1963characterization}; however,  the scattering function is not typically known at the receiver. Often, a flat spectrum in the delay-Doppler domain is assumed, which introduces performance degradation due to the mismatch with respect to the true scattering function \cite{zemen2012adaptive}. 

\revised{
In this work, we adopt a V2V channel model in the delay-Doppler domain, using the geometry-based stochastic channel model proposed in \cite{karedal2009geometry}. Our characterization of this model reveals the special structure of the V2V channel components in the delay-Doppler domain. We show that the delay-Doppler representation of the channel exhibits three key regions; within these regions, the channel is a mixture of specular reflections and diffuse components.} While the specular contributions appear all over the delay-Doppler plane \emph{sparsely}, the diffuse contributions are concentrated in specific regions of the delay-Doppler plane.  We see that the channel measurements from a real data experiment also confirm to our analysis of the V2V channel structures in the delay-Doppler domain. 

  In our prior work \cite{michelusi2012uwb,michelusi2012uwb2}, a Hybrid Sparse/Diffuse (HSD) model was presented for a mixture of sparse and Gaussian diffuse components for a static channel, which we have adapted  to estimate a V2V channel \cite{beygi2014geometry}. This approach requires information about the V2V channel such as the statistics, and power delay profile (PDP) of the diffuse and sparse components \cite{michelusi2012uwb}.  Another approach for time-varying frequency-selective channel estimation is via compressed sensing (CS) or sparse approximation based on an $l_1$-norm regularization \cite{bajwa2010compressed, taubock2010compressive, carbonelli2007sparse}. 
These algorithms perform well for channels with a small number of scatterers or, clusters of scatterers. For V2V channels, diffuse contributions  from reflections along the roadside will degrade the performance of CS methods that only consider element-wise sparsity\cite{taubock2010compressive, zemen2012adaptive}. 

Herein, we exploit the particular structure of the V2V channel, inspired by recent work in 2D sparse signal estimation \cite{ stojnic2009reconstruction, eldar2009robust}, to design a novel joint element- and group-wise sparsity estimator to estimate the 2D time-varying V2V channel using received data. Our proposed method provides a general machinery to solve the joint sparse structured estimation problem with a broad class of regularizers that promote sparsity. We show that our proposed algorithm covers both well-known convex and non-convex regularizers such as smoothly clipped absolute deviation (SCAD) regularizers \cite{fan2001variable},  and the minimax concave penalty (MCP) \cite{zhang2010nearly} that were proposed for element-wise sparsity estimation. We also present a general way to design a proper regularizer function for joint sparsity problem in our previous work \cite{GlobecomPaper}. Our method can be applied to scenarios beyond V2V channels.

Recent algorithms for hierarchical sparsity (sparse groups with sparsity within the groups) \cite{Sprechmann2011CHiLasso, chartrand2013nonconvex, friedman2010note} also consider a mixture of penalty functions (group-wise and element-wise).  Of particular note is \cite{chartrand2013nonconvex} where a similar nested solution is determined, also in combination with the alternating direction method of multipliers (ADMM) as we do herein. Our modeling assumptions can be viewed as a generalization of their assumptions which results in the need for different methods for proving the optimization of the nested structures.  In particular, \cite{chartrand2013nonconvex} examines a particular proximity operator for which the original regularizing function is never specified\footnote{Note given a particular regularization function, there is a unique proximity operator, but for a given proximity operator there may exist more than one regularizer function.}.  This proximity operator is built by generalizing the structure of the proximity operator for the $l_p$ norm.  In contrast,  we begin with a general class of regularizing functions, we specify the properties needed for such functions (allowing for both convex and non-convex functions). Thus, our proof methods rely only on the properties induced by these assumptions.  Furthermore, our results are also applicable to the problem of hierarchical sparsity \cite{Sprechmann2011CHiLasso, chartrand2013nonconvex, friedman2010note}. We observe that the results in \cite{Sprechmann2011CHiLasso, chartrand2013nonconvex, friedman2010note} cannot be applied to non-convex regularizer functions such as  SCAD and MCP due to their concavity and a non-linear dependence on the regularization parameter.    

To find the optimal solution of our joint sparsity objective function,  we take advantage of the alternating direction method of multipliers  \cite{eckstein2012augmented}, which is a very flexible and efficient tool for optimization problems whose objective functions are the combination of multiple terms. Furthermore,  we use the proximity operator \cite{combettes2011proximal} to show that the estimation can be done by using simple thresholding operations in each iteration, resulting in low complexity.

We also address the channel leakage effect due to finite block length and bandwidth for channel estimation in the delay-Doppler plane. In \cite{taubock2010compressive}, the basis expansion for the scattering function is optimized to compensate for the leakage which can degrade performance. The resulting expansion in \cite{taubock2010compressive} is computationally expensive. Herein, we take an alternative view and show that with the proper sampling resolution in time and frequency, we can explicitly derive the leakage pattern and robustify the channel estimator with this knowledge at the receiver to improve the sparsity, compensate for leakage, and maintain modest algorithm complexity. Our overall approach can lead to a performance gain of up to $10$~dB over previously proposed algorithms.

The main contributions of this work are as follows:
\begin{enumerate}
\item A general framework for joint sparsity estimation problem is proposed, which covers a broad class of regularizers including convex and non-convex functions. Furthermore, we show that the solution for joint sparse estimation problem is computed by applying the element-wise and group-wise structure in a nested fashion using simple thresholding operations.

\item We provide a simple model for the V2V channel in the delay-Doppler plane, using the geometry-based stochastic channel modeling proposed in \cite{karedal2009geometry}. We characterize the three key regions in the delay-Doppler domain with respect to the presence of sparse specular and diffuse components. This structure is verified by experimental channel measurement data, as presented in Section \ref{ExperimentSection}.

\item The leakage pattern is explicitly computed and a compensation procedure proposed.

\item A low complexity joint element- and group-wise sparsity V2V channel estimation algorithm is proposed exploiting the aforementioned channel model and optimization result.

\item We use extensive numerical simulation and experimental channel measurement data to investigate the performance of the proposed joint sparse channel estimators and show that our method outperforms classical and compressed sensing methods \cite{zemen2012adaptive, beygi2014geometry, bajwa2010compressed, taubock2010compressive}
\end{enumerate}

The rest of this paper is organized as follows. In Section II, we review some definitions from variational analysis and present our key optimization result for a joint sparse and group sparse signal estimation. In Section III,  the  system model for V2V communications is presented. In Section IV, the geometry-based V2V channel model is developed. The observation model and leakage effect are computed in Section V. In Section VI, the channel estimation algorithm for the time-varying V2V channel model using joint sparsity structure is presented. In Section VII, we provide simulation results and compare the performance of the estimators. In Section VIII, the real channel measurements are provided to confirm the validity of the channel model and numerical simulation. Finally, Section IX concludes the paper. We present our proofs, region specifier algorithm, and analysis of our proximal ADMM in the Appendices.

\emph{Notation}: We denote a scalar by $x$, a column vector by $\boldsymbol{x}$, and its $i$-th element with ${x}[i]$. Similarly, we denote a matrix by $\bf{X}$ and its $(i, j)$-th element by ${X}[{i,j}]$. The transpose of $\bf{X}$ is given by ${\bf X}^{T}$ and its conjugate transpose by ${\bf X}^{H}$. A diagonal matrix with elements $\bf{x}$ is written as $\text{diag}\{\bf{x}\}$ and the identity matrix as $\bf{I}$. The set of real numbers by $\mathbb{R}$, and the set of complex numbers by $\mathbb{C}$. The element-wise (Schur) product is denoted by $\odot$.

\section{Jointly Sparse Signal Estimation:\\ Optimization Result} \label{OptimizationSec}
In this section, we propose a unified framework using proximity operators \cite{rockafellar1998variational} to solve the optimization problem imposed by a structured sparse\footnote{A jointly sparse signal in this paper, is a signal that has both element-wise and group-wise sparsity.} signal estimation problem.  Then, we apply this machinery to estimate the V2V channel, exploiting the group- and element-wise sparsity structure discovered in Section \ref{Ushape}.

\revised{Proximal methods have drawn increasing attention in the signal processing (e.g., \cite{combettes2011proximal}, and numerous references therein) and the machine learning communities (e.g., \cite{bach2012optimization, yu2013decomposing}, and references therein),  due to their convergence rates (optimal for the class of first-order techniques) and their ability to accommodate large, non-smooth, convex (and non-convex) problems. In proximal algorithms, the base operation is evaluating the proximal operator of a function, which involves solving a small optimization problem. These sub-problems can be solved with standard methods, but they often admit closed-form solutions or can be solved efficiently with specialized numerical methods. Our main theoretical contribution is stated in Theorem \ref{MainThm} in Section II.B. In this theorem, we show that, for our proposed class of regularization functions, the nested joint sparse structure can be recovered by applying the element-wise sparsity and group-wise sparsity structures  in a nested fashion, see Eq. \eqref{mainnestedeq}.}

\subsection{Proximity Operator}
We start with the definition of a \emph{proximity operator} from variational analysis \cite{rockafellar1998variational}. 
\begin{defn}\label{ProOperate} 
Let $\phi({\bf a};\lambda)$ be a continuous real-valued function of ${\bf a} \in \mathbb{R}^N$, the
\emph{proximity operator} $P_{\lambda, \phi}({\bf b})$ is defined as 
\begin{align}
\label{Neginbala} 
P_{\lambda, \phi}({\bf b}) :=& \underset{{\bf a} \in \mathbb{R}^N}{\operatorname{argmin}}
\left\{\frac{1}{2}\|{\bf b} - {\bf a}\|_2^2 + \phi({\bf a}; \lambda )\right\},
\end{align}
where  ${\bf b} \in \mathbb{R}^N$ and $\lambda >0$.
\end{defn}
\begin{rem}\label{rem22}
If $\phi(.)$ is a separable function, \emph{i.e.}, $\phi({\bf a};\lambda)
= \sum_{i=1}^{N}f\left(a[i];\lambda\right)$. Then, 
$ [P_{\lambda, \phi}({\bf a})]_i = P_{\lambda,f} (a[i])$.
\end{rem}

\begin{rem} \label{rem23} If the objective function $J({\bf a}) = \frac{1}{2}\|{\bf b} - {\bf a}\|_2^2
+  \phi({\bf a};\lambda)$ is a strictly convex function, the proximity operator of $\phi({\bf a};
\lambda)$  admits a unique solution.
\end{rem}

\begin{rem} \label{rem24} Furthermore, $P_{\lambda,
{\rm \phi}}({\bf b}) $ is characterized by the inclusion 
\begin{equation}
\forall ({\bf a}^*,{\bf b}),\,\,\,{\bf a}^{*}=  P_{\lambda, {\rm \phi}}({\bf b})
 \iff {\bf a}^*-{\bf b}\in \partial \phi({\bf a}^*;\lambda),
\end{equation}
where $\partial \phi(.)$ is the sub-gradient of the function $\phi$ \cite{rockafellar1998variational}.
\end{rem}
 Note that $\phi$ does not need
to be a convex or differentiable function to satisfy the conditions noted in Remarks \ref{rem23} and \ref{rem24}.  Proximity operators have a very natural interpretation in terms of denoising \cite{combettes2011proximal, rockafellar1998variational}. Consider the problem of estimating a vector $\mat{a}\in \mathbb{R}^N$ from an observation $\mat{b}\in \mathbb{R}^N$, namely $\mat{b} = \mat{a} + \mat{n}$
 where $\mat{n}$ is additive white Gaussian noise. If we consider the regularization function $\phi({\bf a}; \lambda )$ as the \emph{prior} information about the vector $\mat{a}$, then $P_{\lambda, \phi}({\bf b}) $ can be interpret as a \emph{maximum a posteriori} (MAP) estimate of the vector $\mat{a}$ \cite{gribonval2011should}.   
\revised{Two well-known types of prior information about the structure of  a vector are element-wise sparsity and group-wise sparsity structures.  A $N$-vector $\mathbf{a}$ is element-wise sparse, if the number of non-zero  (or larger than some threshold) entries in the vector is small compared to the length of the vector. \newline
To define the group-wise sparsity, let us consider $\left\{\mathcal{I}_i\right\}_{i=1}^{N_g}$ be a partition of the index set $\{1, 2, \dots, N\}$ to $N_g$ groups, namely $\cup_{i=1}^{N_g} \mathcal{I}_i =  \{1, 2, \dots, N\}$ and $\mathcal{I}_i \cap \mathcal{I}_j =\varnothing $ for $\forall i\neq j$. We define  group vectors as follows 
\begin{align}
{\bf a}_i[k]=\begin{cases}{\bf a}[k] & k\in \mathcal{I}_i\\ 0& k\notin \mathcal{I}_i\end{cases}
\end{align}
for $i=1, \dots, N_g$, where $N_g$ is the total number of group vectors. Based on above definition, we have ${\bf a} = \sum_{i=1}^{N_g}{\bf a}_i$ and the nonzero elements of ${\bf a}_i$ and ${\bf a}_j$ are non-overlapping for $i \neq j$.  The vector ${\bf a}$ is called a group-wise sparse vector, if the number of group vectors ${\bf a}_i$ for $i=1, \dots, N_g$ with non-zero $l_2$-norm (or $l_2$-norm larger than some threshold) is small compared to the total number of group vectors,  $N_g$.}

\subsection{Optimality of the Nesting of Proximity Operators}\label{SUbsection111}
We consider the estimation of the vector $\mat{a}$ from vector $\mat{b}$ as noted above. Furthermore, suppose that  the vector $\mat{a}$ is a jointly sparse vector. The desired optimization problem is as follows
\begin{align}\label{KeyOpti}
\hat{\bf a} = \underset{{\bf a}\in \mathbb{R}^N}{\operatorname{argmin}} \left\{\frac{1}{2}\left\| {\bf{b}}-{\bf{a}}\right\|_2^2+\phi_g({\bf{a}}, \lambda_g )+ \phi_e({\bf{a}}; \lambda_e)\right\},
\end{align}
where  $\phi_g({\bf{a}}; \lambda_g)$ is a regularization term to induce group sparsity and $\phi_e({\bf{a}}; \lambda_e) $ is a term to induce the element-wise sparsity. In general, the weighting parameters, $\lambda_g > 0, \lambda_e >0$, can be selected from a given range via cross-validation, by varying one of the parameters and keeping the other fixed \cite{friedman2010note}. 
We further consider penalty functions $\phi_g({\bf{a}};\lambda_g) $ and $\phi_e({\bf{a}};\lambda_e)$ of the form:
\begin{align}
\phi_g({\bf{a}}; \lambda_g) &= \sum_{j=1}^{N_g}f_g\left(\|{\bf{a}}_{j}\|_2;\lambda_g\right),\, \text{and}\\ \phi_e({\bf{a}}; \lambda_e) &= \sum_{i=1}^{N}f_e\left({{a}}[i];\lambda_e\right),
\end{align}
where  $f_g:\mathbb{R}\rightarrow \mathbb{R}$ and $f_e:\mathbb{R}\rightarrow \mathbb{R}$ are continuous functions to promote sparsity on groups and elements, respectively, $N$ is the length of  vector ${\bf a}$, and $N_g$ is the number of groups in vector ${\bf a}$.  Our goal here is to derive the solution of the optimization problem in \eqref{KeyOpti} using the proximity operators of the functions $f_e$ and $f_g$. 
We state the conditions imposed on the regularizers, in terms of the univariate functions $f_g(x;\lambda)$ and $f_e(x; \lambda)$ to promote sparsity and also to control the stability of the solution of the optimization problem in \eqref{KeyOpti}.
\begin{enumerate}[]
\item {\bf Assumption I:} For $k \in \{e,g\}$
\begin{enumerate}[i.]
\item $f_k$ is a non-decreasing function of $x$ for $x \geq 0$; $f_k(0;\lambda)=0$; and $f_k(x;0) =0$.
\item $f_k$ is differentiable except at $x=0$.
\item For $\forall z \in \partial f_k(0;\lambda)$, then $|z|\le \lambda$, where   $\partial f_k(0;\lambda)$ is the subgradient\footnote{A precise definition of subgradients of functions is given in \cite{rockafellar1998variational}, page 301.} of $f_k$ at zero.
\item There exists a $\mu\le \frac{1}{2}$ such that the function $f_k(x;\lambda)+\mu x^2$,  is convex. 
\item  $f_g$ is a \emph{homogeneous} function, \emph{i.e.}, $f_g(\alpha x; \alpha \lambda)
= \alpha^2 f_g(x;\lambda)$ for $\forall \alpha>0$.
\item $f_e$ is a \emph{scale invariant} function, \emph{i.e.}, $f_e(\alpha x;\lambda)
= f_e(x;\alpha\lambda) = \alpha f_e(x;\lambda)$ for $\forall \alpha>0$.
\end{enumerate}
\end{enumerate}
It can be observed that conditions (i), (iv), (v), and (vi) ensure the existence of the minimizer of the optimization problem in Eq. \eqref{KeyOpti}, and they induce norm properties on the regularizer function. Assumption (ii) promotes sparsity,  \revised{(iii) controls the stability of the solution in Eq. \eqref{KeyOpti} and guarantees the optimality of solution of optimization problem in Eq. (4) or ($P_0$) in Section VI}. Finally Assumption (iv) enables the inclusion of many non-convex functions in the optimization problem.   Note that  the scale invariant property of $f_e$ implies that $f_e$ also satisfies (v) (is a homogeneous function). 

Many pairs of regularizer functions satisfy Assumption I. For instance, the  $l_1$-norm, namely $f_g(x;\lambda_g) = |x|$ and $f_e(x;\lambda_e) = \lambda_e|x|$, satisfy Assumption I (see Appendix \ref{RegularizerSectionVer}). We note that two recently popularized non-convex functions,  SCAD and MCP regularizers, also satisfy Assumption I.  It is worth pointing out that SCAD and MCP are more effective in promoting sparsity than the $l_p$ norms. The SCAD regularizer \cite{fan2001variable}, is given by
\begin{align}\label{SCADRegula}f_g(x;\lambda)=\begin{cases}\lambda |x| & \text{for}\,\, |x|\le \lambda \\-\frac{x^2-2{\mu_S}\lambda |x|+\lambda^2}{2({\mu_S}-1)} & \text{for}\,\, \lambda<|x|\le {\mu_S}\lambda \\\frac{({\mu_S}+1)\lambda^2 }{2} &\text{for}\,\,|x|>{\mu_S}\lambda\end{cases},\end{align}
where ${\mu_S}> 2$ is a fixed parameter, and the MCP regularizer \cite{zhang2010nearly}, is
\begin{align}\label{MCPRegula}
f_g(x;\lambda)={\rm sign}(x) \int_{0}^{|x|}\left(\lambda-\frac{z}{{\mu_M}} \right)_+\, dz,
\end{align}
where $(x)_+ = \max(0,x)$ and ${\mu_M}> 0$ is a fixed parameter.  In Appendix \ref{RegularizerSectionVer}, we show that Assumption I is met for SCAD and MCP.
The following theorem is our main technical result that presents the solution of the optimization problem in \eqref{KeyOpti} based on the proximity operators of the functions $f_g$ and $f_e$.

\begin{thm}\label{MainThm}Consider functions $f_e$ and $f_g$ that satisfy Assumption I. The optimization problem in \eqref{KeyOpti}, can be decoupled as follows
  \begin{eqnarray}
 {\hat{\bf{a}}}_i =  \underset{{\bf a}_i \in \mathbb{R}^{N}}{\operatorname{argmin}} 
 \left\{\frac{1}{2}\left\| {\bf{b}}_i-{\bf{a}}_i\right\|_2^2+g({\bf a}_i;\lambda_{g})+E({\bf a}_i;\lambda_{e})\right\}, 
\end{eqnarray} 
where the index $i=1\dots N_g$ denotes the group number,
$ g({\bf a}_i;\lambda_{g})  =  f_g \left({\|{\bf a}_i \|_2} ; \lambda_{g} \right),$ and
$E({\bf a}_i;\lambda_{e}) =  \sum_{j}  f_e \left( {{\bf a}_i[j]}; \lambda_{e} \right).$
Then,
\begin{align}\label{mainnestedeq}
\hat{\bf{a}}_i =  \prox{\lambda_{g}, g} \left( \prox{\lambda_{e}, E}({\bf b}_i) \right),
\end{align}
where $ \prox{\lambda_{g}, g}$ and $\prox{\lambda_{e}, E}$ are the proximity operators of $E$ and $g$, respectively (see Definition\,\ref{ProOperate}) and can be written as
\begin{align}
P_{\lambda_{g}, g}\left({\bf b}\right) &=\frac{P_{\lambda_{g}, {f_g}}\left(\|{\bf b}\|_2\right)}{\|{\bf b}\|_2} {\bf b},\\ [P_{\lambda_{e}, E}\left( {\bf b}\right)]_j &= P_{\lambda_{e}, f_e}\left( {b[j]}\right).
\end{align}
\end{thm}
The proof is provided in Appendix \ref{The3}. This result states that within a group, joint sparsity is achieved by first applying the element proximity operator and then the group proximity operator. We observe the $f_g$ and $f_e$ can be chosen structurally different, the resultant complexity is modest, and we can use our result with any optimization algorithm based on proximity operators, \emph{e.g.}, ADMM \cite{eckstein2012augmented}, proximal gradient methods \cite{rockafellar1998variational}, proximal splitting methods \cite{combettes2011proximal}, and so on.

\subsection{Proximity Operator of Sparse-Inducing Regularizers}
In this section, we compute closed form expressions for the proximity
operators of the sparsity inducing regularizers introduced in Section \ref{SUbsection111},
 using Definition \ref{ProOperate} and Remarks \ref{rem22} to \ref{rem24}. All of the aforementioned regularizers
satisfy the condition noted  in Remark \ref{rem23} due to  property
(iv) in Assumption I.  Using Definition \ref{ProOperate}, we can compute
the proximity operators of the $l_p$, SCAD, and MCP regularizers.
The proximity operator for the $l_p$-norm is given by, $P_{\lambda,
f_g}({x}) = {\rm sign}(x) \max\{0,\lambda u\}$, where $u^{p-1}+\frac{u}{p}=\frac{|x|}{\lambda
p}$. For $p=1$, \emph{i.e.}, $f_e(x;\lambda) = \lambda |x|$, the resulting operator is often called {\em soft-thresholding} (see {\em e.g.} \cite{friedman2010note}),
\begin{align}
P_{\lambda, f_e}(x)={\rm sign}(x) \max\{0, |x|-\lambda\}.
\end{align}

\begin{figure}
\centering
\includegraphics[width = 3 in]{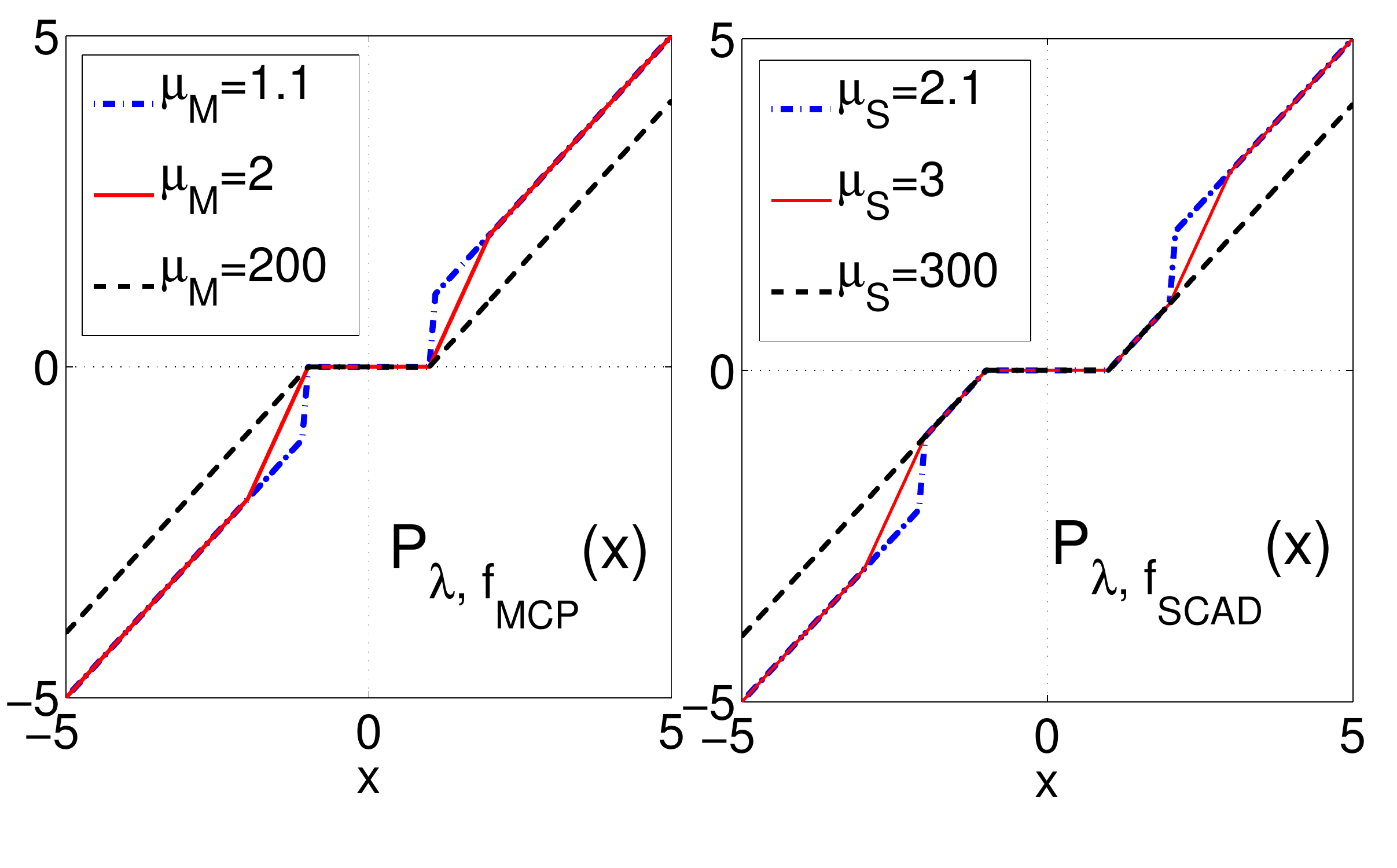}
\vspace{-0.2in}
\caption{MCP and SCAD proximity operators. Here $\lambda=1$ is considered.}
\vspace{-0.1in}
\label{FigProximal}
\end{figure}
The closed form solution of the proximity operator for the SCAD regularizer \cite{fan2001variable}
can be written as,
\begin{align}\label{multiSCAD}
P_{\lambda, f_g}(x)=\begin{cases}0, & \text{if}\,\, |x|\le \lambda\\
x-{\rm sign}(x)\lambda,& \text{if}\,\, \lambda \le |x|\le 2\lambda\\
\frac{x-{\rm sign}(x)\frac{{\mu_S}\lambda}{{\mu_S}-1}}{1-\frac{1}{{\mu_S}-1}} & \text{if}\,\,
2\lambda<|x|\le {\mu_S}\lambda \\ x &\text{if}\,\,|x|>{\mu_S}\lambda\end{cases},
\end{align}
 and finally the proximity operator for the MCP regularizer \cite{zhang2010nearly} is 
\begin{align}\label{multiMCP}
P_{\lambda, f_g}(x)=\begin{cases}0, & \text{if}\,\, |x|\le \lambda \\
\frac{x-{\rm sign}(x){\lambda}}{1-\frac{1}{{\mu_M}}} & \text{if}\,\, \lambda<|x|\le
{\mu_M}\lambda \\ x &\text{if}\,\,|x|>{\mu_M}\lambda\end{cases}.
\end{align}
In Fig.~\ref{FigProximal}, MCP and SCAD proximity operators are depicted
for $\lambda=1$ and three different values of the parameters ${\mu_S}$ and ${\mu_M}$.
It is clear that when ${\mu_S}$ and ${\mu_M}$ are large, both SCAD and MCP operators
behave like the soft-thresholding operator (for $x$ smaller than ${\mu_S}\lambda$
and ${\mu_M}\lambda$, respectively).

In the sequel, first we model a V2V communication system. Then, we show that V2V channel representation in the delay-Doppler domain has both element-wise and group-wise sparsity structures.
Finally, we apply our key optimization result derived in this section to estimate the V2V channel using  an ADMM algorithm.


\section{Communication System Model}\label{SystemModel}

We will consider communication between two vehicles as shown in Fig.~\ref{fig:V2V_geometry}. 
The transmitted signal $s(t)$ is generated by the modulation of the transmitted pilot sequence $s[n]$ onto the transmit pulse $p_{t}(t)$ as,
\begin{align}\label{Eq11}
s(t) = \sum_{n=-\infty}^{+\infty}s[n]p_{t}(t-nT_{s}),
\end{align}
where $T_{s}$ is the sampling period. Note that this signal model is quite general, and encompasses OFDM signals as well as single-carrier signals. The signal $s(t)$ is transmitted over a linear, time-varying, V2V channel. The received signal $y(t) $ can be written as,
\begin{align}\label{Eq12}
y(t) = \int_{-\infty}^{+\infty}h\left(t,\tau\right)s(t-\tau)\,d\tau + z(t).
\end{align}
Here, $h(t,\tau)$ is the channel's time-varying impulse response, and $z(t)$ is a complex white Gaussian noise. At the receiver, $y(t)$ is converted into a discrete-time signal using an anti-aliasing filter $p_{r}(t)$. That is,
\begin{align}\label{Eq13}
y[n] = \int_{-\infty}^{+\infty}y(t)p_{r}(nT_{s}-t)\, dt.
\end{align}
The relationship between the discrete-time signal $s[n]$ and received signal $y[n]$, using Eqs.~\eqref{Eq11}--\eqref{Eq13}, can be written as,
\begin{align}\label{Eq14}
y[n] = \sum_{m = -\infty}^{+\infty}h_{l}\left[n,m\right]s[n-m] + z[n],
\end{align}
where $h_{l}[n,m]$ \footnote{\revised{The subscript $``l"$ hereafter denotes the channel with leakage. We discuss the channel leakage effect with more details in Section \ref{se:LekageComp}.}} is the discrete time-delay representation of the observed channel, which is related to the continuous-time channel impulse response $h(t,\tau)$ as follows,
\begin{equation*}
h_{l}[n,m] = \iint\limits_{-\infty}^{+\infty} h\left(t+nT_{s},\tau\right)p_{t}(t-\tau+mT_{s})p_{r}(-t)\,dt d\tau.
\end{equation*}
With some loss of generality, we assume that $p_r(t)$ has a root-Nyquist spectrum with respect to the sample duration $T_s$, which implies that $z[n]$ is a sequence of \emph{i.i.d} circularly symmetric complex Gaussian random variables with variance $\sigma_z^2$, and that $h_{l}[n,m]$ is causal with maximum delay $M-1$, \emph{i.e.}, $h_{l}[n,m]=0$ for $m \ge M $ and $m < 0$. 
We can then write 
\begin{align}
\label{ReceivedRxDDD} y[n] 
=\sum_{m = 0}^{M-1}\left( \sum_{k=-K}^{K}  {H}_{l}\left[k,m\right]e^{j\frac{2\pi nk}{2K+1}}\right) s[n-m]  + z[n],& \\   \text{for} \quad n= 0, 1, . . ., N_r-1& \nonumber, 
\end{align}
where $2K+1 \ge N_r$, and
\begin{align}\label{DDDRep}
H_l[k, m] = \frac{1}{2K+1} \sum_{n=0}^{N_r -1} h_l[n, m] e^{-j\frac{2\pi nk}{2K+1}}, \,\text{for}\, |k|\le K,
\end{align}
is the \emph{discrete delay-Doppler, spreading function} of the channel. \revised{Here, $N_r$ denotes the total number of received signal samples used for the channel estimation}.

\section{Joint Sparsity Structure of  V2V Channels}\label{Ushape}
In this section, we adopt the V2V geometry-based stochastic channel model from~\cite{karedal2009geometry} and analyze the structure such a model imposes on the delay-Doppler spreading function. The V2V channel model considers four types of multipath components (MPCs): $(i)$ the effective \emph{line-of-sight} (LOS) component, which may contain the ground reflections, $(ii)$ discrete components generated from reflections of discrete mobile scatterers (MD), \emph{e.g.}, other vehicles, $(iii)$ discrete components reflected from discrete static scatterers (SD) such as bridges, large traffic signs, etc., and $(iv)$ diffuse components (DI).  Thus, the V2V channel impulse response can be written as
 \begin{align}
\label{channelModel} h(t,\tau) =& h_{LOS}(t,\tau)+\sum_{i = 1}^{N_{MD}}h_{MD,i}(t,\tau) + \sum_{i= 1}^{N_{SD}}h_{SD,i}(t,\tau) +\sum_{i = 1}^{N_{DI}}h_{DI,i}(t,\tau),
 \end{align}
where $N_{MD}$ denotes the number of discrete mobile scatterers, $N_{SD}$ is the number of discrete static scatterers and $N_{DI}$ is the number of diffuse scatterers, respectively. 
Typically, $N_{DI}$ is much larger than $N_{SD}$ and $N_{MD}$~\cite{karedal2009geometry}.
In the above representation, the multipath components can be modeled as  
\begin{align}\label{singlepathmodel}
h_{i}(t,\tau) = \eta_{i}\delta(\tau-\tau_{i})e^{-j2\pi \nu_{i}t},
\end{align}
where $\eta_{i}$ is the complex channel gain, $\tau_{i}$ is the delay, and $\nu_{i}$ is the Doppler shift associated with path $i$ and $\delta(t)$ is the Dirac delta function.
\begin{figure}
  \centering
  \includegraphics[scale=0.3]{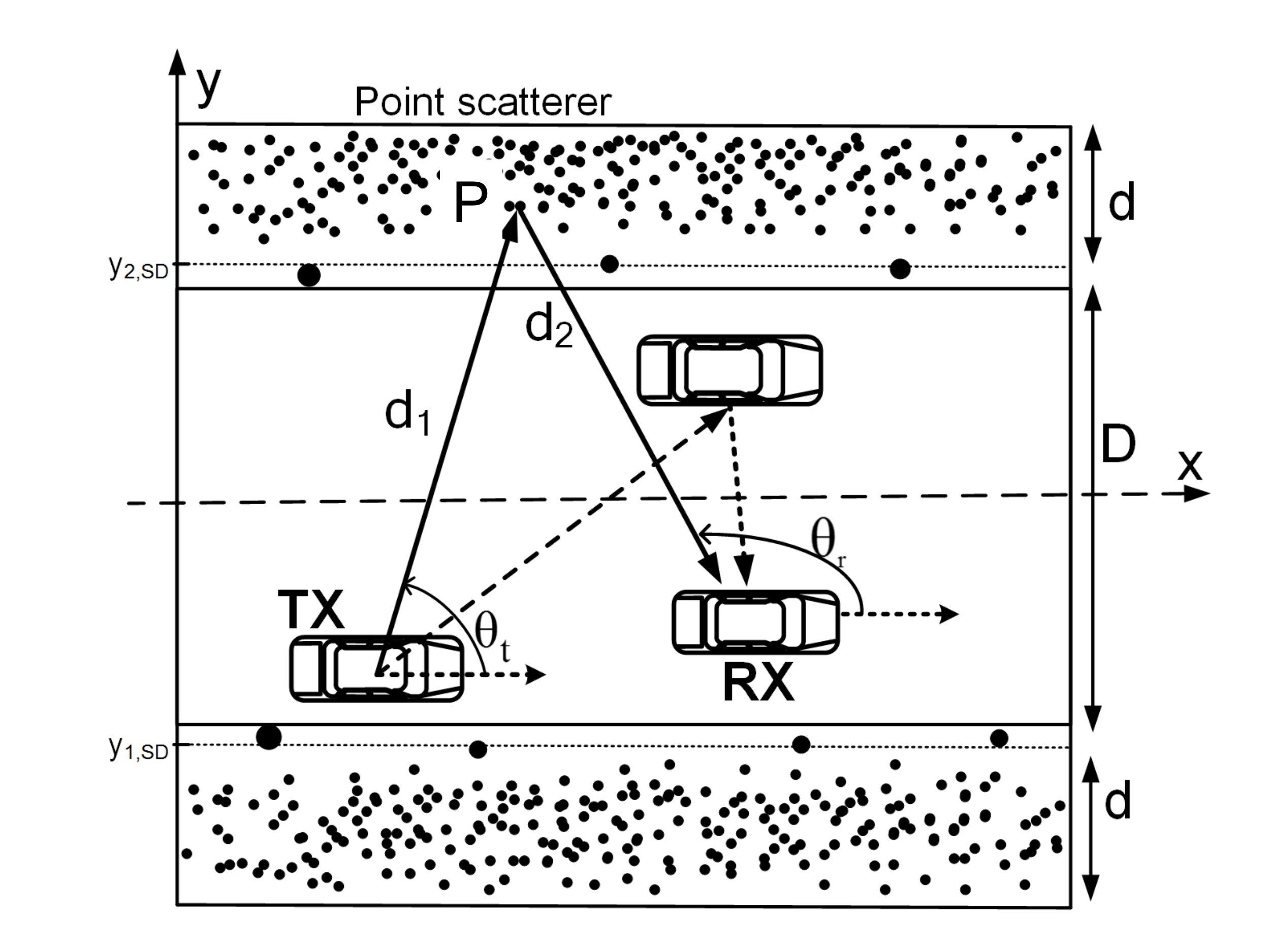}
  \caption{Geometric representation of the V2V channel. The shaded areas on each side of the road contain static discrete (SD) and diffuse (DI) scatters, while the road area contains both SD and moving discrete (MD) scatterers.}
  \vspace{-0.15in}
  \label{fig:V2V_geometry}
\end{figure}
\revised{The channel description in \eqref{channelModel} and  \eqref{singlepathmodel} explicitly models distance-dependent pathloss and scatterer parameters [7]. We assume that these parameters can be approximated as time-invariant over the pilot signal duration. This is a reasonable assumption in practical systems as will be illustrated in Figures \ref{PracticalFig0} and \ref{PracticalFig34} for the experimental channel measurement data.}
The spatial distribution of the scatterers and the statistical properties of the complex channel gains are specified in \cite{karedal2009geometry} for rural and highway environments. It is shown that the channel power delay profile is exponential. Further details about the spatial evolution of the gains can be found in \cite{karedal2009geometry, borhani2013correlation}.  
In geometry-based stochastic channel modeling, point scatterers are randomly distributed in a geometry according to a specified distribution. The position and speed of the scatterers, transmitter, and receiver determine the delay-Doppler parameters for each MPC, which in turn, together with the transmitter and receive pulse shapes, determine $H_l[k, m]$. 

We next determine the delay and Doppler contributions of an ensemble of point scatterers of type $(i)$-$(iv)$ for the road geometry depicted in Fig.~\ref{fig:V2V_geometry}. If vehicles are assumed to travel parallel with the $x$-axis, the overall Doppler shift for the path from the transmitter (at position TX) via the point scatterer (at position P) to the receiver (at position RX) can be written as \cite{karedal2009geometry}
\begin{equation}\label{DopplerEq}
\nu \left(\theta_{t}, \theta_{r}\right) = \frac{1}{\lambda_{\nu}}\left[\left(v_{T}-v_{P}\right)\cos \theta_{t} + \left(v_{R}-v_{P}\right)\cos \theta_{r} \right],
\end{equation}
where $\lambda_{\nu}$ is the wavelength, $v_T$, $v_P$, and $v_R$ are the speed of the transmitter, scatterer, and receiver, respectively, and $\theta_t$ and $\theta_r$ is the angle of departure and arrival, respectively. The path delay is
\begin{equation}\label{delayEq}
\tau = \frac{d_{1}+d_{2}}{c_{0}},
\end{equation}
where $c_0$ is the propagation speed, $d_{1}$ is the distance from TX to P, and $d_{2}$ is the distance from P to RX. The path parameters $\theta_t$, $\theta_r$, $d_1$, and $d_2$ are easily computed from TX, P, and RX. The delay and Doppler parameters of each component ($i$)-($iv$) can now be specified by Eqs.~\eqref{DopplerEq} and \eqref{delayEq}.

\paragraph*{LOS}
If it exists, the most significant component of the V2V channel is the line of sight (LOS) component, which will have delay and Doppler parameters
$\tau_{0} = \frac{d_{0}}{c_{0}}$ and  $\nu_{0} = \frac{1}{\lambda_{\nu}}\left(v_{T}-v_{R} \right)\cos(\theta),$
where $d_0$ is the distance from TX and RX and $\theta$ is the angle between the x-axis (\emph{i.e.}, the moving direction) and the line passing through TX and RX. 

\begin{figure}[t]
 \centering
 \includegraphics[width=2.6in]{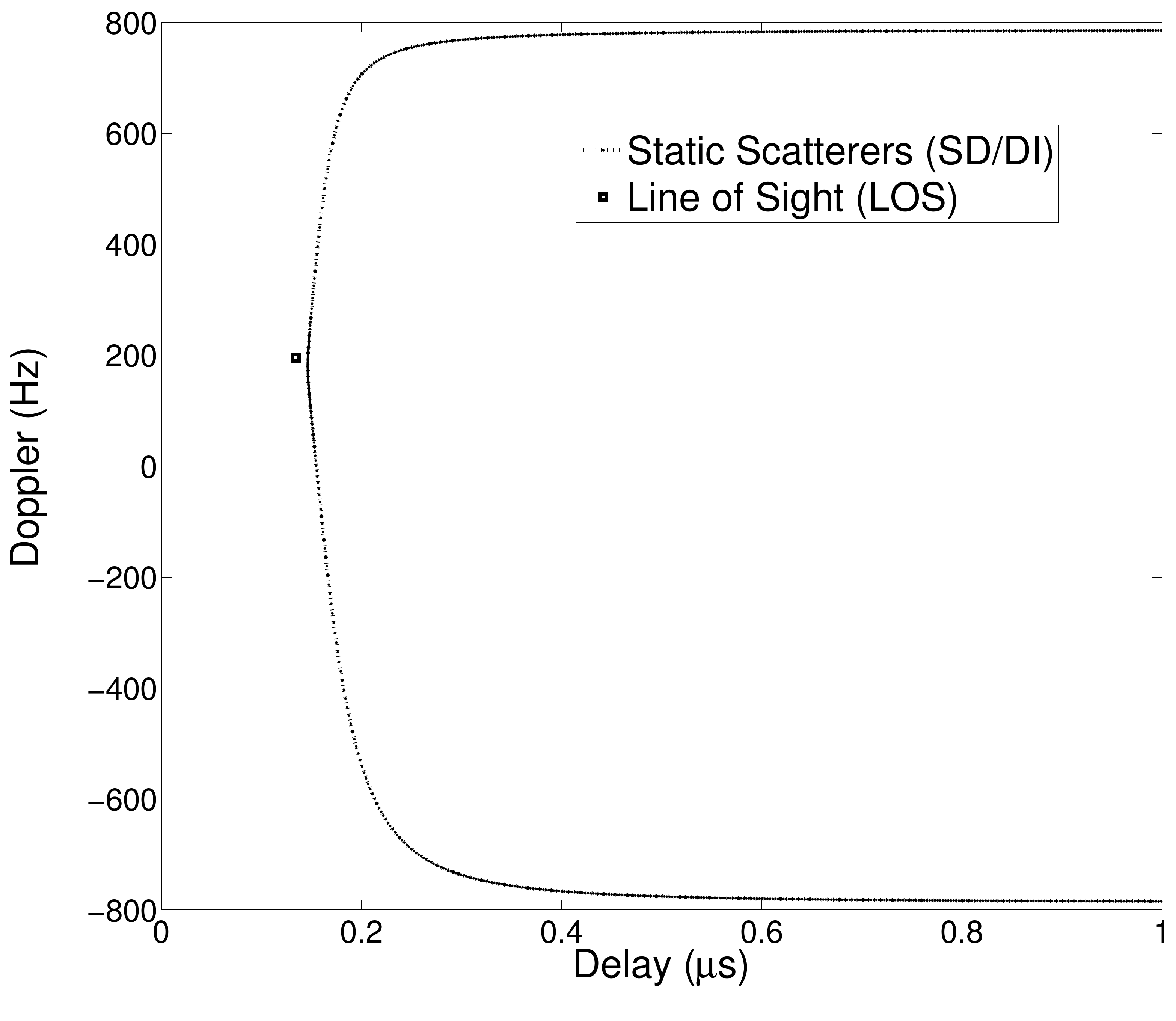}
 \vspace{-0.1in}
 \caption{Delay-Doppler contribution of line-of-sight (LOS) and from static scatterers (SD/DI) placed on a parallel line beside the road.}\vspace{-0.17in}
 \label{fig:SingleUcurve}
\end{figure}

\paragraph*{Diffuse Scatterers} 
The diffuse (DI) scatterers are static ($v_P = 0$) and uniformly distributed in the shadowed regions in Fig.~\ref{fig:V2V_geometry}. Suppose we place a static scatterer at the coordinates $(x, y)$. The delay-Doppler pair, $\left(\tau(x, y), \nu(x, y)\right)$, for the corresponding MPC can be calculated from Eqs.~\eqref{DopplerEq} and~\eqref{delayEq}. If we fix $y=y_0$ and vary $x$ from $-\infty$ to $+\infty$, the delay-Doppler pair will trace out a U-shaped curve in the delay-Doppler plane, as depicted in Fig.~\ref{fig:SingleUcurve}. 

\begin{figure}
  \centering
  \includegraphics[width=2.8in]{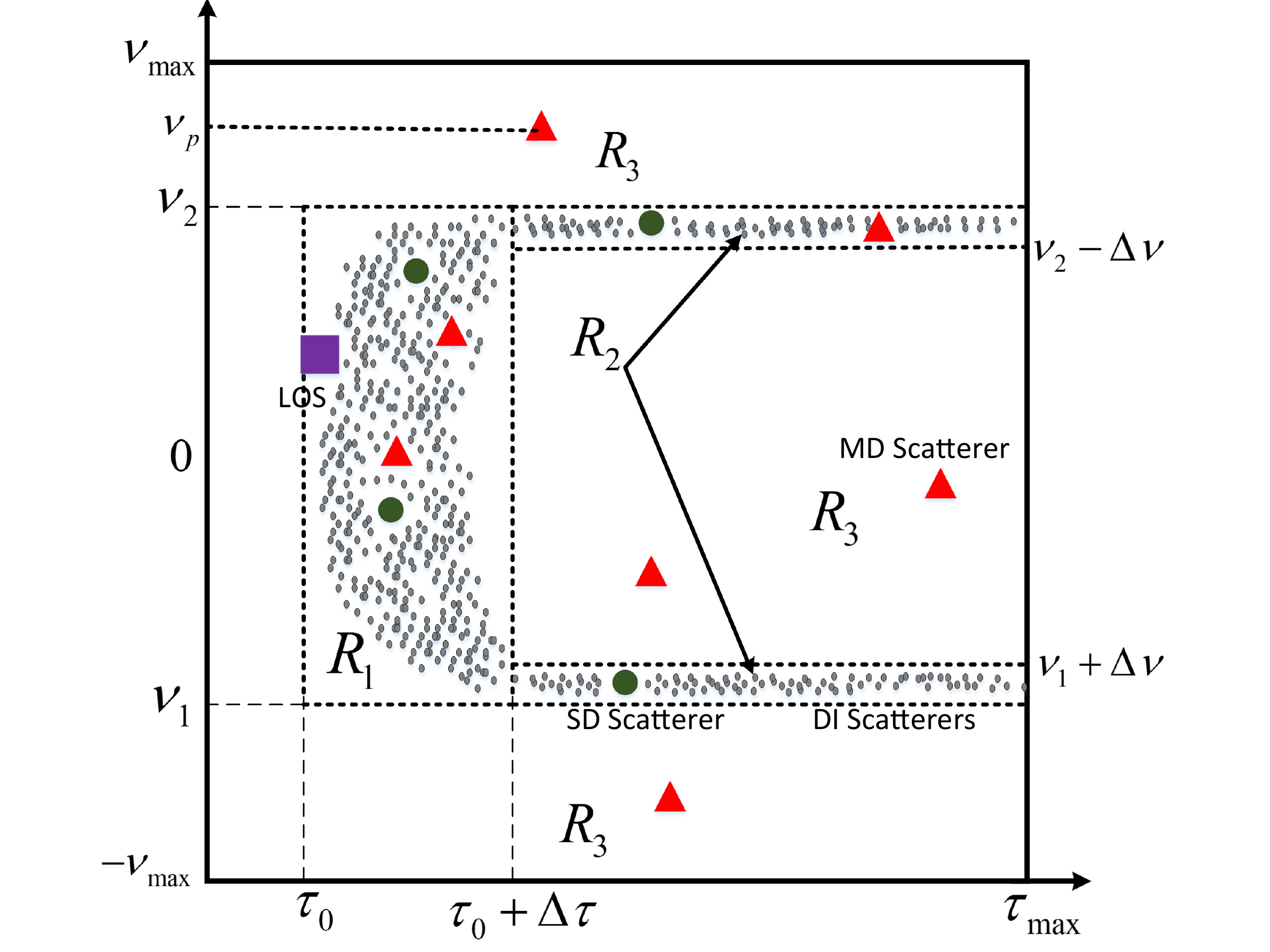}
  \vspace{-0.1in}
  \caption{Delay-Doppler domain representation of V2V channel. Delay-Doppler spreading function for diffuse components is confined to a U-shaped area. }\vspace{-0.15in}
  \label{FigDParts}
\end{figure}

Repeating this procedure for the permissible $y$-coordinates for the DI scatterers, 
$|y_0| \in [D/2, d+D/2]$, will result in a family of curves that are confined to a U-shaped region, see Fig.~\ref{FigDParts}. Hence, the DI scatterers will result in multipath components with delay-Doppler pairs inside this region. The maximum and minimum Doppler values of the region is easily found from Eq.~\eqref{DopplerEq}. In fact, it follows from Eq.~\eqref{DopplerEq}  that the Doppler parameter of an MPC due to a static scatterer will be confined to the symmetric interval $[-\nu_{S}, \nu_{S}]$, where 
 $\nu_{S}=\frac{1}{\lambda_{\nu}} (v_T + v_R).$      

\paragraph*{Static Discrete Scatterers} 
The static discrete (SD) scatterers can appear outside the shadowed regions in Fig.~\ref{fig:V2V_geometry}. In fact, the $y$-coordinates of the SD scatterers are drawn from a Gaussian mixture consisting of two Gaussian pdfs with the same standard deviation $\sigma_{y,SD}$ and means $y_{1,SD}$ and $y_{2,SD}$~\cite{karedal2009geometry}. The delay-Doppler pair for an MPC due to an SD scatterer can therefore appear also outside the U-shaped region in Fig.~\ref{FigDParts}. However, since the SD scatterers are static, the Doppler parameter is in the interval $[-\nu_{S}, \nu_{S}]$, \emph{i.e.}, the same interval as for the diffuse scatterers.

\paragraph*{Mobile Discrete Scatterers}
We assume that no vehicle travels with an absolute speed exceeding $v_{\text{max}}$. It then follows from Eq.~\eqref{DopplerEq} that the Doppler due to a mobile discrete (MD) scatterer is in the interval $[-\nu_{\text{max}}, \nu_{\text{max}}]$, where
 $  \nu_{\text{max}}=\frac{4v_{\text{max}}}{\lambda_{\nu}}.      $
For example, in Fig.~\ref{FigDParts}, the Doppler shift $\nu_{p}$  is due to an MD scatterer (vehicle) that travels in the oncoming lane ($v_P < 0$). 

Based on the analysis above, we can conclude that the delay-Doppler parameters for the multipath components can be divided into three regions,
\begin{align*}
R_1 &\identi \left\{(\tau, \nu)\in \mathbb{R}^2: \tau \in (\tau_0, \tau_0+\Delta \tau), \nu\in(-\nu_{S}, \nu_{S})\right\}\\
R_2 &\identi \left\{(\tau, \nu)\in \mathbb{R}^2:   
          \tau\in[\tau_0+\Delta\tau, \tau_{\text{max}}], 
          |\nu| \in [\nu_{S}-\Delta\nu, \nu_{S}) \right\}\\
R_3 &\identi \left\{(\tau, \nu)\in \mathbb{R}^2:   
          \tau\in[\tau_0, \tau_{\text{max}}], 
          |\nu|\le \nu_{\text{max}}] \right\} \setminus (R_1\cup R_2),
\end{align*}
where $\tau_{\text{max}}-\tau_0$ is the maximum significant excess delay for the V2V channel. Here, $\Delta\tau$ and $\Delta\nu$ are chosen such that the contributions from all diffuse scatterers are confined to $R_1\cup R_2$. 
The exact choice of $\Delta \tau$ is somewhat arbitrary. In this paper, we consider a thresholding rule to compare the noise level and channel components, which results in a particular choice of $\Delta \tau$; the method is specified  in Appendix~\ref{CoarseEstimation}.
However, regardless of method, once $\Delta \tau$ is chosen, we can compute the height $\Delta \nu$ of the two strips that make up  $R_2$. This can be done by placing an ellipsoid with its foci at the transmitter and receiver such that the path from the transmitter to the receiver via any point on the ellipsoid has propagation delay $\tau_0+\Delta \tau$. By computing the associated Doppler along the part of ellipsoid that is in the diffuse region (\emph{i.e.}, in the strips just outside the highway, see Fig.~\ref{fig:V2V_geometry}), we can determine the smallest absolute Doppler value among them as $\nu'$ and calculate $\Delta\nu$ as $\Delta\nu=\nu_{S}-\nu'$.
In Appendix \ref{CoarseEstimation}, we present a data driven approach to approximate $\Delta\nu$. 
Note that the regions gather channel components with similar behavior. Region $R_1$, contains the LOS, ground reflections, and (strong) discrete and diffuse components due scatterers near the transmitter and receiver. In Region $R_2$, the delay-Doppler contribution of static discrete and diffuse scatterers from farther locations appear. Region $R_3$ contains contributions from moving discrete and static discrete scatters only. 
\begin{rem}
In Fig.~\ref{FigDParts}, we see that there are sparse contributions from the SD and MD scatterers in all regions $R_1$, $R_2$, and $R_3$. However, clusters of DI components are confined to $R_1\cup R_2$.  \revised{Therefore, V2V channel components can be divided into two main clusters. One is the element-wise sparse components (mobile and static discrete scatterers) that are distributed in all three regions $R_1$, $R_2$, and $R_3$,  and the other one is the group-wise sparse components (diffuse components) that are located in regions $R_1$ and $R_2$. We note that, in V2V channels due to the geometry of the channel and antenna heights, there exists more diffuse components compared to the cellular communication.  Therefore, a proper V2V channel estimation algorithm should consider the estimation of diffuse components with higher quality compared to the cellular communication channel. Since the diffuse components are located in a specific part of delay-Doppler domain and the rest of the channel support in the delay-Doppler is essentially zero, we partition the channel into group vectors and take advantage of group-wise sparsity of the channel to enhance the accuracy of the estimate of the diffuse components.} 
In Sec.~\ref{ITmeansHere}, we propose a method based on joint element-wise and group-wise sparsity and Theorem \ref{MainThm} of Section \ref{SUbsection111} to estimate the discrete delay-Doppler representation of V2V channel exploiting this structure.
\end{rem}
\section{Observation Model and Leakage Effect}\label{se:LekageComp}
In this section, we show how pulse shaping and a finite-length training sequence can be taken into account when formulating a linear observation model of the V2V channel. We assume that $p_t(t)$ and $p_r(t)$ are causal with support $[0, T_{\text{supp}})$. The contribution to the received signal from one of the $1+N_{MD}+N_{SD}+N_{DI}$ terms in \eqref{channelModel} is of the form
\begin{align}
  s(t) * h_i(t, \tau) = \sum_{l=-\infty}^\infty s[l] \eta_{i}e^{j2\pi\nu_i t} p_t(t-lT_s-\tau_i)\approx \sum_{l=-\infty}^\infty s[l] \eta_{i}e^{j2\pi\nu_i (lT_s + \tau_i)} p_t(t-lT_s-\tau_i),
\end{align}
where the approximation is valid if we make the (reasonable) assumption that $\nu_i T_{\text{supp}} \ll 1$, and $*$ denotes convolution. If we let $p(t) = p_t(t) * p_r(t)$, we can write the contribution after filtering and sampling as
\begin{align} 
  y_i[n] &= \left. s(t) * h_i(t, \tau) * p_r(t)\right|_{t=nT}\\ \nonumber&= \sum_{l=-\infty}^\infty s[l] \eta_{i}e^{j2\pi\nu_i(lT_s +\tau_i)} p((n-l)T_s-\tau_i)\\\nonumber
  &= \sum_{m=-\infty}^\infty s[n-m] \eta_{i}e^{j2\pi\nu_i((n-m)T_s +\tau_i)} p(mT_s-\tau_i) ,
\end{align}
and identify
$ h_i[n, m] = \eta_{i}e^{j2\pi\nu_i((n-m)T_s +\tau_i)} p(mT_s-\tau_i).$
Suppose we have access to $h_i[n,m]$ for $n=0, 1, \ldots, N_r-1$ and
let \revised{$\omega_{2K+1} = \exp(j2\pi/(2K+1))$}. The $(2K+1)$-point DFT of $h_i[n,
m]$, where we choose $(2K+1)\ge N_r$, is
\begin{align}\nonumber
H_{l,i}[k, m] &= \frac{1}{2K+1}\sum_{n=0}^{N_r-1} h_i[n, m]\omega_{2K+1}^{-nk}  \qquad k \in\mathcal{K}\\
&= \eta_{i} e^{-j2\pi\nu_i(mT_s-\tau_i)}  p(mT_s-\tau_i) w(k, \nu_i),
\end{align}
where $\mathcal{K}=\{0, \pm 1, \pm 2, \ldots, \pm K\}$ and $w(k, x)$
is the $(2K+1)$-point DFT of a discrete-time complex exponential with
frequency $x$ and truncated to $N_r$ samples
\begin{align}
w(k, x) = 
\begin{cases}
  \frac{N_r}{2K+1}, \quad x = k/(2K+1) \quad& {} \\
  \frac{e^{-j\pi\left(\frac{k}{2K+1}-x\right) (N_r-1)}}{2K+1} \frac{\sin\left(\pi\left(\frac{k}{2K+1}-x\right) N_r \right) }{\sin\left(\pi\left(\frac{k}{2K+1}-x\right) \right)},
& \text{otherwise}
\end{cases}
\end{align}
We note that the leakage in the delay and Doppler plane is due to the
non-zero support of $p(\cdot)$ and $w(\cdot, \cdot)$. The leakage with respect to Doppler decreases with
the observation length, $N_r$, and the leakage with respect to  delay decreases with
the bandwidth of the transmitted signal. 
\revised{We compute the (effective) channel coefficients at time sample $m_i = \left[\frac{\tau_i}{T_s}\right]$ and Doppler sample $k_i = \left[\nu_i T_s(2K+1) \right]$, where $[.]$ indicates the closest integer number. Note that the true channel parameters $\tau_i$ and $\nu_i$ are not restricted to integer multiples of a sampling interval. However, we do seek to estimate the effective channel after appropriate sampling. Thus, at the receiver side, to compensate (but not perfectly remove) for the channel leakage, the effective channel components at those sampled times are computed to equalize the channel (which may be different from the actual channel components)}. We can then write
\begin{align}
H_{l,i}[k, m] = \eta_{i} g[k, m, k_i, m_i], \qquad k \in \mathcal{K}, m\in\mathcal{M}
\end{align}
where $\mathcal{M} = \{0,1, \ldots, M-1\}$ and
\begin{align}
g[k, m, k', m'] = \omega_{2K+1}^{-k'(m-m')}w(k-k', 0) p((m-m')T_s). 
\end{align}

Due to the linearity of the discrete Fourier transform, we can
conclude that the channel with leakage is given by
\begin{align}\label{LeakageEffect00}
H_{l}[k,m] &= \sum_{i} H_{l,i}[k,m] = \sum_{i} \eta_{i} g[k,m, k_i, m_i], 
\end{align}
where the summation is over the LOS component and all the
$N_{MD}+N_{SD}+N_{DI}$ scatterers. The channel without leakage is 
\begin{align}\label{NoLeakageEffect}
H[k,m] &= \sum_{i} \eta_{i} \delta[k-k_i] \delta[m-m_i], 
\end{align}
where $\delta[n]$ is the Kronecker delta function.
The channels in
\eqref{LeakageEffect00} and \eqref{NoLeakageEffect}  can be
represented for $k\in\mathcal{K}$ and $m \in \mathcal{M}$
by the vectors
${\bf{x}}_{l}\in\mathbb{C}^{N}$ and
${\bf{x}} \in \mathbb{C}^{N}$, respectively, where $N =
|\mathcal{K}||\mathcal{M}|=(2K+1)M$, as
\begin{align}
{\bf{x}}_{l} &= 
\vstack \begin{bmatrix}
  H_l[-K, 0] &  \cdots &  H_l[-K, M-1]\\
  \vdots &  \cdots & \vdots\\
  H_l[K, 0] &   \cdots &   H_l[K, M-1]\\
\end{bmatrix}\\  \label{Vectorizing}
{\bf{x}}&= 
\vstack \begin{bmatrix}
  H[-K, 0] &   &   H[-K, M-1]\\
  \vdots &  \cdots & \vdots\\
  H[K, 0] &  \cdots &   H[K, M-1]\\
\end{bmatrix}.
\end{align}
where $\vstack(\mathbf{H})$ is the vector formed by stacking the columns of ${\bf H}$ on the
top of each other.
The relationship between $\mathbf{x}_l$ and $\mathbf{x}$ can be
written as
\begin{equation}\label{LeakageRel}
{\bf{x}}_{l} = {\bf{G}}{\bf{x}},
\end{equation}
where ${\bf{G}}\in \mathbb{C}^{N\times N}$ is defined as
\begin{align}
{\bf{G}} &=
\begin{bmatrix}
  \vstack(\mathbf{G}_0) & \vstack(\mathbf{G}_1) & \cdots & \vstack(\mathbf{G}_{N-1})
\end{bmatrix}\label{Gmatrix}\\ \nonumber
\mathbf{G}_j &=
\begin{bmatrix}
  g[-K, 0, k', m'] &   \cdots &   g[-K, M-1, k',  m']\\
  \vdots & \cdots & \vdots\\
  g[K, 0, k', m'] &   \cdots &   g[K, M-1, k',  m']\\
\end{bmatrix},
\end{align}
where the one-to-one correspondence between $j=0, 1, \ldots, (2K+1)M-1$
and $(k', m')\in\mathcal{K}\times\mathcal{M}$ is given by $j = m'(2K+1) +
k' + K$. 

The structure of ${\bf{G}}$ is a direct consequence of how we
vectorize ${H}_{l}[k,m]$ in \eqref{Vectorizing}. If we consider
an alternative way of vectorizing ${H}_{l}[k,m]$, then the leakage matrix
${\bf{G}}$ needs to be recomputed accordingly, by appropriate
permutation of the columns and rows of leakage matrix defined
in~\eqref{Gmatrix}. 
As $K$, $M$, and the pulse shape are known, ${\bf{G}}$ is completely determined in \eqref{LeakageRel}. Thus, we can
utilize the relationship in \eqref{LeakageRel} to compensate for
leakage.

Consider that the source vehicle transmits a sequence of $N_{r}+M-1$ pilots, $s[n]$, for $n = -(M-1), -(M-2), ..., N_r-1$, over the channel. We collect the $N_{r}$ received samples in a column vector
\begin{align}{\bf{y}} = \left[ y[0], y[1], . . ., y[N_{r}-1]\right]^\intercal.\end{align} Using \eqref{ReceivedRxDDD}, we have the following matrix representation:
\begin{equation}
{\bf{y}} = {\bf{S}}{\bf{x}}_{l}+{\bf{z}},
\label{eq:y_Sxl_z}
\end{equation}
where ${\bf{z}} \sim \mathcal{CN}({\bf{0}}, \sigma_z^2 {\bf{I}}_{N_r})$ is a Gaussian noise vector, and ${\bf{S}}$ is a $ N_{r}\times N $ block data matrix of the form 
\begin{align}
{\bf{S}}= \left[{{\bf{S}}_0, . . . , {\bf{S}}_{M-1}}\right],
\end{align}
where each block ${\bf{S}}_{m} \in \mathbb{C}^{N_{r}\times (2K+1)}$ is of the form 
\begin{align}{\bf{S}} _{m}= \text{diag}\left\{ s[-m], . . . , s[N_{r}-m-1]\right\}{\bf{\Omega}},\end{align}  for $m = 0, 1, \ldots, M-1$, and ${\bf{\Omega}} \in \mathbb{C}^{N_{r}\times (2K+1)}$ is a Vandermonde matrix, ${\bf{\Omega}}[{i,j}]=\omega_{2K+1}^{i(j-K)}$, where 
$i=0, 1\ldots, N_r-1$  and $j=0, 1, \ldots, 2K$.
Finally, by combining~\eqref{eq:y_Sxl_z} and \eqref{LeakageRel} we have   
\begin{equation}
{\bf{y}} = {\bf{S}}{\bf{x}}_{l}+{\bf{z}} = {\bf{A}}{\bf{x}}+{\bf{z}},
\end{equation}
where ${\bf{A}} = {\bf{S}}{\bf{G}}$ and \revised{${\bf{A}}\in \mathbb{C}^{N_{r}\times N}$}.

\section{Channel Estimation}\label{ITmeansHere}
Based on our analysis in Section \ref{Ushape},  the components in the vector ${\bf x}$ have both element- and group-wise sparsity.  
\revised{
Given estimates of the parameters that define regions $R_1$, $R_2$, and $R_3$ (see Section IV and Appendix C), we illustrate how to partition the elements in the channel vector ${\bf x}$ to enforce the group-wise sparsity structure.  Our partitioning is based on the sparsity structure of regions $R_1$, $R_2$, and $R_3$,  and Eq. \eqref{Vectorizing}, which maps the entries of  ${\bf x}$ into channel components, $H[k,m]$ for $k=-K, \dots, K$ and $m=0, \dots, M-1$.
 \begin{figure}[!t]
\centering
    \includegraphics[width=0.45\textwidth]{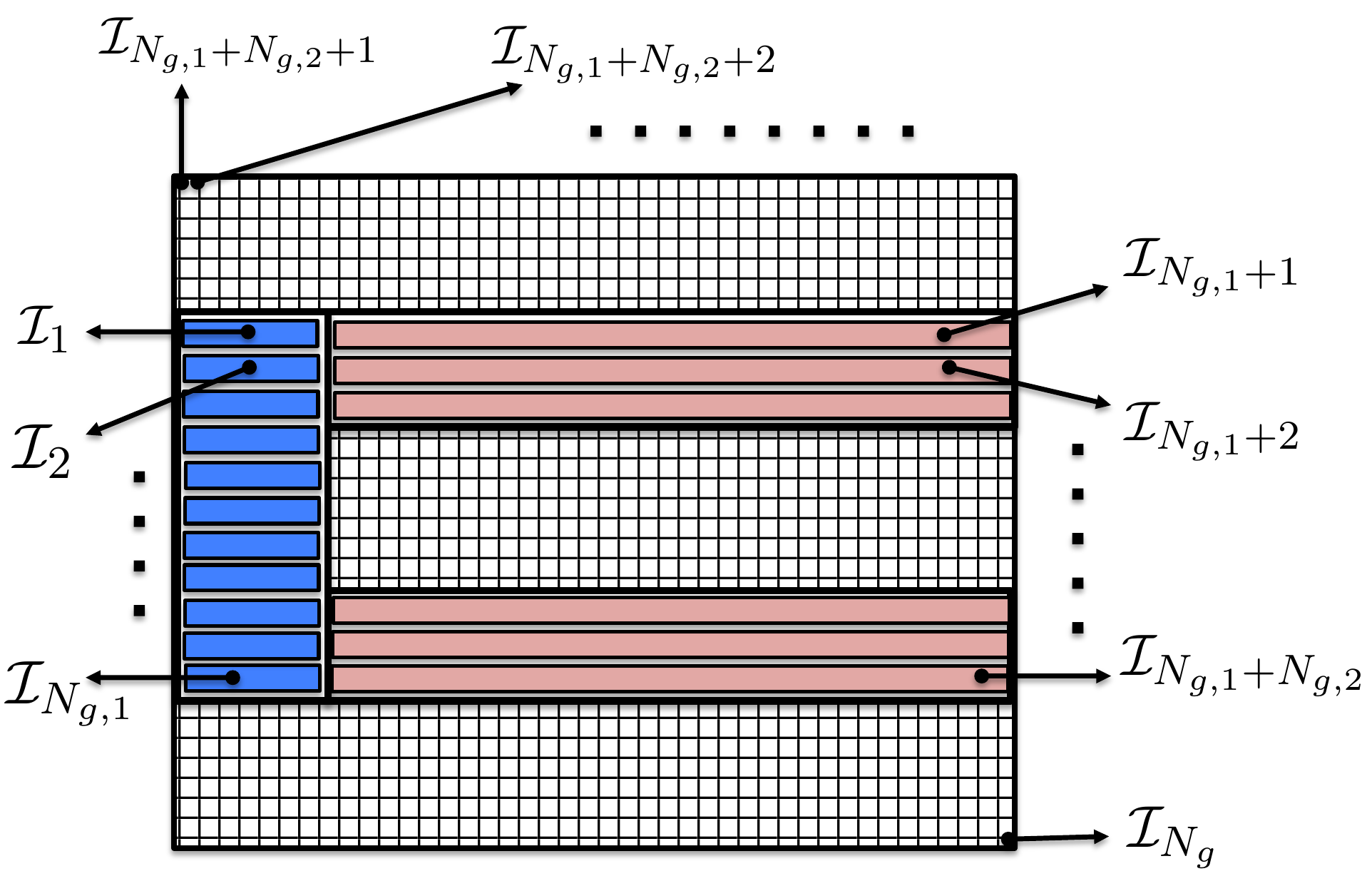}
    \caption{Schematic of the V2V channel vector partitioning for group vectors.}
    \label{fig:groupings}
\end{figure}
To partition the elements in regions $R_1$ and $R_2$, we collect channel components  with a common Doppler value into a single group, \emph{e.g.},  $\mathcal{I}_1, \mathcal{I}_2$, \dots, $\mathcal{I}_{N_{g,1}}$ (blue segments),  and $\mathcal{I}_{N_{g,1}+1} , \dots, \mathcal{I}_{N_{g,1}+N_{g,2}}$  (red segments) as depicted in Fig.\,\ref{fig:groupings}.   For $R_3$, we know that this region contains only the element-wise sparse components, thus we consider each element of ${\bf x}$ in this region as a single  partition, \emph{e.g.},  $\mathcal{I}_{N_{N_g,1} + N_{N_g,2} + 1}, \dots,  \mathcal{I}_{N_g}$, where $N_g = N_{N_g,1} + N_{N_g,2} + N_{N_g,3}$ as depicted in Fig.\,5. Now that we have a partition of all the elements in ${\bf x}$, we can easily determine the group vectors ${\bf x}_i$ as follows:
\begin{align}
{\bf x}_i[k]=\begin{cases}{\bf x}[k] & k\in \mathcal{I}_i\\ 0& k\notin \mathcal{I}_i\end{cases}, \,\,\text{for}\,\, i=1, \dots, N_g
\end{align}
 We know that most of components in $R_3$ are zero (or close to zero) and there is a significant numbers of non-zero diffuse components in region $R_1$ and $R_2$.  Therefore, if we enforce the group sparsity regularization over this partitioning of elements in ${\bf x}$, it will improve the the quality of estimation of diffuse components.}

We next specify the regularizations to exploit the jointly sparse structure of the V2V channel as follows 
\begin{align*}
\hat{\bf x} = \underset{{\bf x} \in \mathbb{C}^N}{\operatorname{argmin}} \left\{\frac{1}{2}\left\| {\bf{y}}-\mat{A}{\bf{x}}\right\|_2^2+\phi_g(|{\bf{x}}|;\lambda_g)+ \phi_e(|{\bf{x}}|;\lambda_e)\right\}, \quad (P_0)
\end{align*}
where  $|\mat{x}| = \left[ |x[1]|, \dots, |x[N]| \right]^T$ with $N = M(2K+1)$ and $|x[i]| = \sqrt{\text{Re}(x[i])^2+\text{Im}(x[i])^2 }$. Here the regularization functions are
\begin{align}
\phi_g(|{\bf{x}}|; \lambda_g) &= \sum_{j=1}^{N_g}f_g\left(\|{\bf{x}}_{j}\|_2;\lambda_g\right),\\ \phi_e(|{\bf{x}}|; \lambda_e) &= \sum_{i=1}^{N}f_e\left(|{{x}}[i]|;\lambda_e\right).
\end{align}

We develop a proximal alternating direction method of multipliers (ADMM) algorithm to solve problem $P_0$. Problem $P_0$ can be rewritten using an auxiliary variable ${\bf{w}}$ as follows,
\begin{align}\label{OptimizationMain}
\nonumber \underset{{\bf x, w} \in \mathbb{C}^N}{\operatorname{min}} \frac{1}{2}\left\| {\bf{y}}-{\bf{A}}{\bf{x}}\right\|_2^2&+\phi_g(|{\bf{w}}|;\lambda_g)+ \phi_e(|{\bf{w}}|;\lambda_e)\\
\text{s.t.}& \, \,{\bf{w}} ={\bf{x}}
\end{align}
For the optimization problem in \eqref{OptimizationMain}, ADMM consists of the following iterations,

\begin{itemize}
\item {\bf Initialize}: $\rho\neq 0$, $\lambda_{\rho g}=\frac{\lambda_g}{\rho}$,  $\lambda_{\rho e}=\frac{\lambda_e}{\rho}$, ${\boldsymbol{\theta}}^{0} ={\bf{w}}^{0} = {\bf 0}$, ${\bf{A}}_0 = \left({\rho^2}{\bf{I}}+{\bf{A}}^{H}{\bf{A}}\right)^{-1}$, and ${\bf{x}}_0 =  {\bf{A}}_0{\bf{A}}^H{\bf{y}}$.  
\item {\bf Update}-{\bf x}:
\begin{align}
{\bf{x}}^{n+1} = {\rho^2}{\bf{A}}_0 \left({\bf{w}}^{n} -{\boldsymbol \theta}_{\rho}^{n}\right) + {\bf{x}}_0,
\end{align}
where ${\boldsymbol \theta}_{\rho}^{n} = \frac{{\boldsymbol \theta}^{n}}{\rho^2}$.
\item {\bf Update}-{\bf w}:  for $ i=1, 2, \dots, N_g,$
\begin{align}\nonumber
{\bf{w}}^{n+1}_i =  \underset{{\bf w}_i}{\operatorname{argmin}} \,\,\frac{1}{2}\left\|{\bf{x}}_i^{n+1}+{\boldsymbol \theta}_{\rho i}^{n}- {\bf w}_i\right\|_2^2\\\label{wstepnon}+ g(|{\bf w}_i|;\lambda_{\rho g})
 +E(|{\bf w}_i|;\lambda_{\rho e})  , \,\,\, 
\end{align}
where  the index $i$ denotes the group number, and
\begin{align}\label{TermSeper001}
E(|{\bf w}_i|;\lambda_{\rho e}) &= 
\sum_{j}f_e\left(\frac{|{w}_i[j]|}{\rho};\lambda_{\rho e}\right),\\ 
\label{TermSeper002}
g(|{\bf w}_i|;\lambda_{\rho g}) &= f_g\left(\frac{\|{\bf{w}}_i\|_2}{\rho};\lambda_{\rho g}\right).
\end{align}

\item {{\bf Update-dual variable}-${\boldsymbol{\theta}}$}:
\begin{align}
  {\boldsymbol \theta}_{\rho}^{n+1} =  {\boldsymbol \theta}_{\rho}^{n} +  \left({\bf{x}}^{n+1}-{\bf{w}}^{n+1}\right).
\end{align}
\end{itemize} 

Details of this derivation are provided in Appendix \ref{gabay1976dualAPEN}. Note that both ${\bf{A}}_0$ and ${\bf{x}}_0$ are known and can be computed in advance. In the update-${\bf{w}}$ step, the index $i$ denotes the group number, thus, this step can be done in parallel for all groups, simultaneously. 

Since the vectors in optimization problem for updating $\mat{w}$ in \eqref{wstepnon} are complex vectors, Theorem \ref{MainThm} in Section \ref{SUbsection111}, cannot directly be applied to find a closed-form solution for this optimization problem. However, in order to apply Theorem \ref{MainThm}, we introduce the following notation and lemma. The vector $\mathbf{w} \in \mathbb{C}^n$
can be written as $\mathbf{w} = |\mathbf{w}|\odot {\rm Phase}(\mathbf{w})$,
where the $n$th element of ${\rm Phase}(\mathbf{w})$ is
$\exp(j{\rm Ang}(\mathbf{w}[n]))$, and 
${\rm Ang}(\mathbf{w}[n])$ is the angle of $\mathbf{w}[n]$ in polar form, \emph{i.e.}, $\mathbf{w}[n]=|\mathbf{w}[n]|\exp(j{\rm Ang}(\mathbf{w}[n]))$, and $\odot $ is component-wise multiplication (Schur product).

\begin{lem}\label{LemComplexCOnv}
For any $\mathbf{c}\in \mathbb{C}^N$
\begin{align}
\amin{\mathbf{z}\in \mathbb{C}^N} \|\mathbf{c} - \mathbf{w}\|_2^2 = 
{\rm Phase}(\mathbf{c})\odot\amin{\mathbf{|w|}\in \mathbb{R}^N} \| \,|\mathbf{c}| - |\mathbf{w}|\,\|_2^2 .
\end{align}
\end{lem}
The proof is provided in Appendix \ref{ProofLemma}. Since the two last terms in~\eqref{wstepnon} are independent of the
phase of $\mathbf{w}_i$, we use Lemma~\ref{LemComplexCOnv} to write the
$i$th group problem in~\eqref{wstepnon} as
 \begin{align}\label{Equation34}
{\bf{w}}^{n+1}_i ={\rm Phase}\left({\bf{x}}_i^{n+1}+{\boldsymbol \theta}_{\rho i}^{n}\right) \odot\left(\underset{|{\mat{w}_i}|}{\operatorname{argmin}} \frac{1}{2}\left\||{\bf{x}}_i^{n+1}+{\boldsymbol \theta}_{\rho i}^{n}|- |{\mat{w}_i}|\right\|_2^2+ g(|{\bf w}_i|;\lambda_{\rho g}) 
+E(|{\bf w}_i|;\lambda_{\rho e}) \vphantom{\underset{|{\mat{w}_i}|}{\operatorname{argmin}}}\right)
\end{align} 
Now, the vectors in the optimization problem in \eqref{Equation34} involve only real vectors, therefore, the solution for this optimization problem can be directly computed using Theorem \ref{MainThm}. 
We determine a closed form solution for the update-${\bf w}$ step in Corollary \ref{MainCoro}, below, using the proximity operators of the univariate functions $f_e$ and $f_g$. This update rule is a direct consequence of Theorem~\ref{MainThm}.

\begin{coro}\label{MainCoro}The second step, update-${\bf w}$, can be performed as follows
\begin{align}
{\bf{w}}^{n+1}_{i}= P_{\lambda_{\rho g}, g}\left( P_{\lambda_{\rho e}, E}\left(\left|{\bf{x}}_i^{n+1}+{\boldsymbol \theta}_{\rho i}^{n}\right| \right)\right)\
\odot {{\rm Phase}}\left({\bf{x}}_i^{n+1}+{\boldsymbol \theta}_{\rho i}^{n}\right),
\end{align}
where  $E$ and $g$ are defined in Equations \eqref{TermSeper001} and \eqref{TermSeper002}.
\end{coro}
  
Based on Corollary \ref{MainCoro}, the update-${\bf w}$ step only depends on the proximity operators of the regularizer functions $f_e$ and $f_g$. The proposed algorithm to estimate the channel vector ${\bf x}$ from the received data vector ${\bf y}$ is summarized in Table I.
\begin{table}[!t] \label{Table1}
{\normalsize 
\begin{center}
\caption{Proposed V2V Channel Estimation Method}
\begin{tabular}{|l|}
\hline
{\bf Input:} {\bf{y}}, ${\bf A}$, $\lambda_g$, $\lambda_e$, $\rho$, $n_{\max}$, $\epsilon$.\\
 {\bf Initialize:} $ {\bf{w}}^{0} = {\boldsymbol \theta}_{\rho}^{0} = {\bf 0}$.\\
 {\bf Pre-computation:} ${\bf A}_0 = \left({\rho^2}{\bf{I}}+{\bf{A}}^{H}{\bf{A}}\right)^{-1}$, \\~~~~~~~~~~~~~~~~~~~~~~~${\bf x}^0 = {\bf A}_0 {\bf{A}}^H{\bf{y}}$.\\
{\bf For} $n~ = ~0~ \text{to} ~ n_{\max}-1$\\
~~~~~~~$ {\bf{x}}^{n+1} = {\rho^2}{\bf A}_0\left({\bf{w}}^{n} -{\boldsymbol \theta}_{\rho}^{n}\right)+{\bf x}^0 $\\
~~~~~~~$
{\bf{w}}^{n+1}_{i}= P_{\lambda_{\rho g}, g}\left( P_{\lambda_{\rho e}, E}\left(\left|{\bf{x}}_i^{n+1}+{\boldsymbol \theta}_{\rho i}^{n}\right| \right)\right)$\\
~~~~~~~~~~~~~~~~~~~~~~~~~~~~~~~~~~~$\odot {\rm Phase}\left({\bf{x}}_i^{n+1}+{\boldsymbol \theta}_{\rho i}^{n}\right)
$\\
\quad \quad \quad \quad \quad \quad ~~~~~~~~~~~ for $ i=1, 2, \dots, N_g$
 \\
~~~~~~~${\boldsymbol \theta}_{\rho}^{n+1} =  {\boldsymbol \theta}_{\rho}^{n} +  \left({\bf{x}}^{n+1}-{\bf{w}}^{n+1}\right)$ \\
~~~~~~~{\bf if} $\|{\bf{x}}^{n+1}-{\bf{x}}^{n}\|_2< \epsilon$ then {\bf break}\\
 {\bf End}\\
 {\bf Output:} Vector {\bf{x}}.
  \\
\hline
\end{tabular}
\end{center}}
\end{table}
\revised{The complexity of our joint sparse signal estimation algorithm is as follows: 
{ Step 1)} update ${\bf x}$, incurs a computational complexity of $O(N^2)$, where $N = M(2K+1)$, due to a matrix/vector multiplication.
{Step 2)} updating ${\bf{w}}$, requires $N_g$ group-wise threshold comparisons and $N$ element-wise threshold comparisons. Thus, it has complexity of $O(N+N_g) \approx O(N)$.  Updating the dual variable in Step 3) has $O(N)$ complexity.

Therefore, the overall complexity of our proposed algorithm is $O(N^2)$.   In our algorithm (similar to the other methods), we  compute  a Least-Squares (LS) solution for the initialization. Computation of a least squares (LS) solution can be implemented with complexity $O(N^3)$.  Algorithms such as the Wiener filter, the Hybrid Sparse/Diffuse (HSD) estimator that are designed based on statistical knowledge of channel parameters require a large number of samples to estimate the required covariance matrices. If the correlation matrices are known, then the Wiener filter and the Hybrid Sparse/Diffuse (HSD) estimator have computational complexity in order of $O(N^3)$. }


\section{Numerical Results}\label{SimulationResult}
In this section, we demonstrate the performance gains that can be achieved with our proposed, structured, estimation using both convex and non-convex sparsity-inducing regularizers, in comparison to prior methods such as Wiener filtering \cite{zemen2012adaptive}, the Hybrid Sparse/Diffuse (HSD) estimator \cite{michelusi2012uwb, michelusi2012uwb2, beygi2014geometry}, and the compressed sensing (CS) method \cite{bajwa2010compressed, taubock2010compressive}. 

To simulate the channel, we consider a geometry with  length of $1$ km around the transmitter-receiver pair, road width $D = 50$\,m, and the width of the diffuse strip around the road $d=25$\,m, as Fig.~\ref{fig:V2V_geometry}. The locations of the transmitter and receiver are chosen in this geometry with distance $100$\,m to $200$\,m from each other. The speeds of the transmit and receive vehicles are chosen randomly from the interval $[60,160]$ (km/h), the speed limits for a highway. It is assumed that the number of MD scatterers $N_{MD} = 10$, and their speeds are also chosen randomly from the interval $[60,160]$ (km/h);  we have $N_{SD}=10$ and $N_{DI}=400$,  SD and DI scatterers, respectively \cite{karedal2009geometry}. Using these parameters, we compute the delay and Doppler values for each scatterer. The statistical parameter values for different scatterers are selected as in Table I of \cite{karedal2009geometry}, which are determined from experimental measurements. The scatterer amplitudes were randomly drawn from zero mean, complex Gaussian distributions with three different power delay profiles for the LOS and mobile discrete (MD) scatterers, static discrete scatterers, and diffuse scatterers. We assume that the mean power of the static scatterers is $10$ dB less than the mean power of the LOS and MD scatterers, and the mean power of the diffuse scatterers also is $20$ dB less than the mean power of the LOS and MD scatterers \cite{karedal2009geometry}. Furthermore, we consider $f_c = 5.8$ GHz, $T_s = 10$ ns,  $N_r=1024$, $K=512$, and $M=256$. The pilot samples are drawn from a zero-mean, unit variance Gaussian distribution. The interpolation/anti-aliasing filters $p_t(t) =p_r(t)$ are root-raised-cosine filters with roll-off factor $0.25$ and $T_{\rm supp}=1 \mu$s. The required regularization parameters were found by trial and error using a cross-validation algorithm on the data \cite{picard1984cross}.  Performance is measured by the normalized mean square error (NMSE), normalized by the mean energy of the channel coefficients. The \mbox{NMSE} is defined as $E\left\{{\|\hat{{\bf x}}-{\bf x}\|_2^2}/{\|{\bf x}\|_2^2}\right\}$, where $\hat{{\bf x}}$ is the estimated channel vector and SNR defined as $\mbox{SNR}={E\{\|{\bf y}-{\bf z}\|_2^2\}}/{E\{\|{\bf z}\|_2^2\}}$. \\
 
\begin{figure}
\centering
\includegraphics[width = 2.8in]{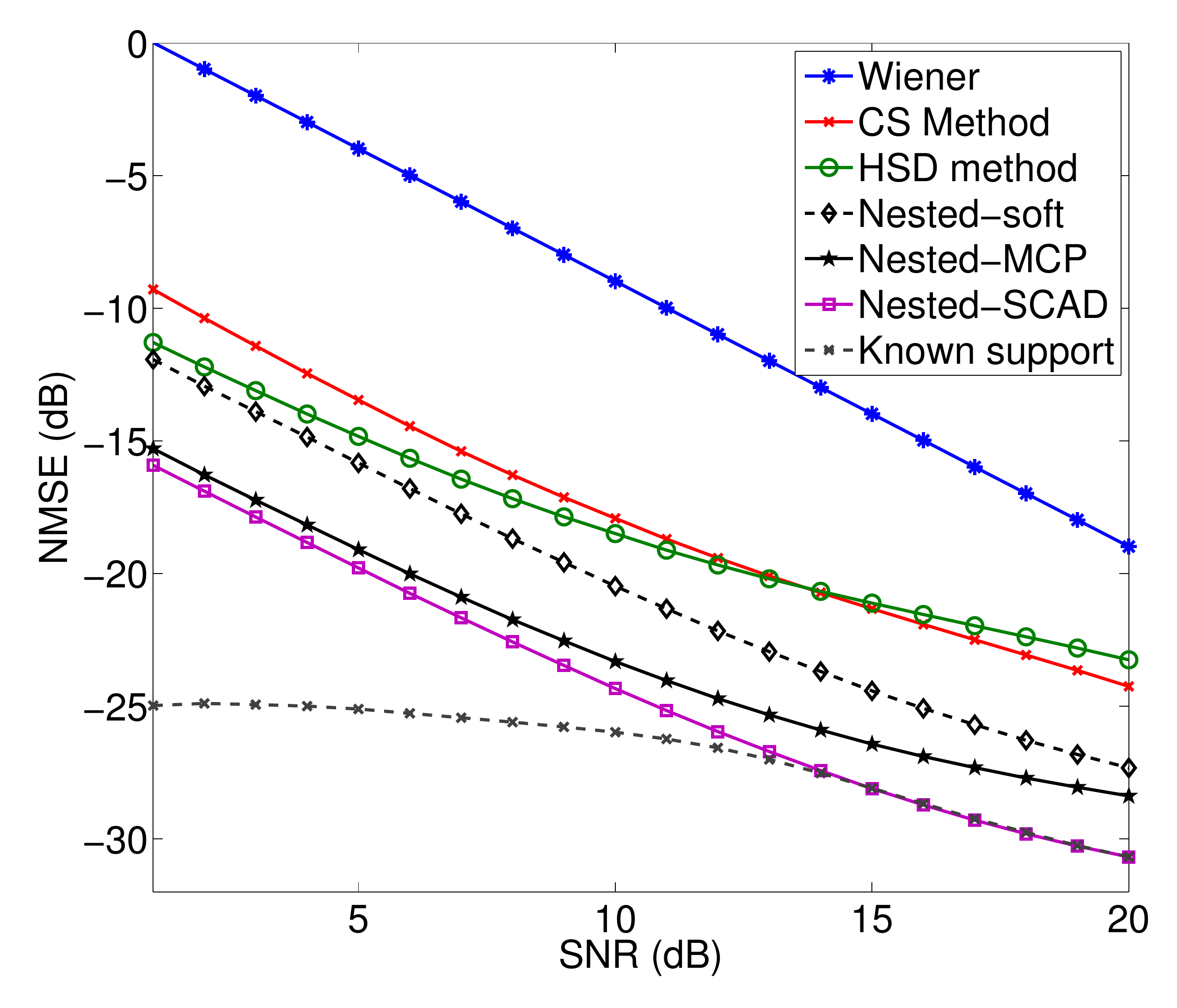}
\vspace{-0.1in}
\caption{Comparison of NMSE v.s SNR for Wiener filter estimator  \cite{zemen2012adaptive}, HSD estimator \cite{beygi2014geometry, michelusi2012uwb}, CS Method \cite{bajwa2010compressed, taubock2010compressive}, and proposed method with different regularizes, \emph{i.e.}, Nested-soft \cite{Sprechmann2011CHiLasso, chartrand2013nonconvex}, Nested-MCP, and Nested-SCAD regularizers.}\vspace{-0.1in} 
\label{FigDParts21}
\end{figure}

Fig.~\ref{FigDParts21} depicts the MSE using our proposed estimator, our previously proposed hybrid sparse/diffuse (HSD) estimator as adapted to V2V channels\cite{beygi2014geometry},  a compressed sensing (CS) method \cite{bajwa2010compressed, taubock2010compressive}, and the Wiener based estimator of \cite{zemen2012adaptive}. {We have also included a curve (known support) corresponding to the case that the support (location of non-zero components) of vector ${\bf x}$ is known when we apply our joint sparse estimation method.} This curve provides a lower bound for the performance of our proposed estimator.   
The HSD estimator considers the channel components as a summation of sparse and diffuse components, {\emph{i.e.}} ${\bf x}={\bf x}_s+{\bf x}_d$. Sparse components, ${\bf x}_s$, are modeled by the element-wise product of an unknown deterministic amplitude, ${\bf a}_s$, and a random Bernoulli vector, ${\bf b}_s$,  ${\bf x}_s={\bf a}_s \odot{\bf b}_s$. Furthermore, the diffuse components are assumed to follow a Gaussian distribution with exponential profile, ${\bf x}_d \sim \mathcal{N}({\bf 0}, {\boldsymbol{\Sigma_d}})$ where ${\boldsymbol{\Sigma_d}}$ is diagonal and is the Kronecker product of covariance  matrices of the channel component vectors with common Doppler values \cite{beygi2014geometry}. The profile parameters are retrieved using the expectation-maximization algorithm \cite{michelusi2012uwb}. 

The HSD model estimation procedure can be stated as: first, the location of sparse components, \emph{i.e.}, the Bernoulli vector, are determined as $\hat{b}_{s}[k]=\mathds{1}\left(|{{x}}_{LS}[k]|^{2}\ge \gamma \left( {({\Sigma}_e[{k,k}])}^{-1}+{{ \Sigma}_d[{k,k}]}\right) \right)$, where $\mathds{1}(.)$ is the indicator function, ${\bf x}_{LS} = \left({\rho^2}{\bf{I}}+{\bf{A}}^{H}{\bf{A}}\right)^{-1}{\bf{A}}^H{\bf{y}}$ is the regularized LS estimation, ${\bf \Sigma}_e$ is the covariance matrix of the noise vector after LS estimation and $\gamma$ is a known parameter \cite{beygi2014geometry}. Then, the sparse components are computed as $\hat{{\bf{x}}}_{s} = {\bf{x}}_{LS}\odot \hat{\bf b}_{s},$ and finally the diffuse components can be estimated as
$\hat{\bf{x}}_{d} = {\bf \Sigma}_{d} \left({\bf \Sigma}_{d}+{\bf \Sigma}_e^{-1}\right)^{-1}\left({\bf{x}}_{LS}-\hat{\bf{x}}_{s}\right)$ \cite{beygi2014geometry, michelusi2012uwb}.

 The Wiener based estimator \cite{zemen2012adaptive} estimates the channel as 
 $\hat{\bf x} = {\bf R}_x{\bf A}^H\left({\bf A}{\bf R}_x{\bf A}^H+\sigma_z^2 {\bf I} \right)^{-1}{\bf y}$ with ${\bf R}_x=\mathbb{E}\{{\bf x}{\bf x}^H\}$, which is not known at the receiver side, but is approximated by assuming a delay-Doppler scattering function prototype with flat spectrum in a 2D region as in \cite{bello1963characterization}. The maximum path delay and Doppler in the support of scattering function are considered as $\tau_{\max}=1.5 \,\mu$s and $\nu_{\max}=860$ Hz, respectively. As we will see, this idealized scattering function assumption results in degraded performance. 
  
It is clear from Fig.~\ref{FigDParts21} that there is performance improvement when we consider the joint structural information of the channel in the delay and Doppler domain. For our proposed method, we have considered different types of regularizers. In Fig. ~\ref{FigDParts21}, \emph{Nested-Soft} corresponds to our proposed structured estimator with a soft-thresholding regularizer, \emph{i.e.}, $f_g(x;\lambda_g)  = \lambda_g|x|$ and $f_e(x;\lambda_e) = \lambda_e|x|$ {(Note that this special case of our algorithm is the case considered in\cite{Sprechmann2011CHiLasso} and \cite{chartrand2013nonconvex} for unknown parameter $p=1$)}; \emph{Nested-SCAD} corresponds to the case where $f_g(x;\lambda_g)$ is the SCAD regularizer function with $\mu_S=3$ in Equation \eqref{SCADRegula} and $f_e(x;\lambda_e) = \lambda_e|x|$; and  \emph{Nested-MCP} corresponds to the case where $f_g(x;\lambda_g)$ is the MCP regularizer function with $\mu_M=2$ in Equation \eqref{MCPRegula} and $f_e(x;\lambda_e) = \lambda_e|x|$, respectively. 
The penalty parameters for the simulations have been considered as $\lambda_g/\lambda_e \approx 10$ and $\lambda_g\in [0, 10]$. 
 Fig.~\ref{FigDParts21} shows that the non-convex regularizers improve estimation quality by about $5$ dB  at low SNR and $7$ dB at high SNR values with the same computational complexity compared to the convex soft-thresholding regularizer.   There is also a significant improvement in effective SNR due to the exploitation of the V2V channel structure in the delay-Doppler domain. For instance, to achieve \mbox{MSE} =  $-20$ dB,  the performance curve related to the the structured estimator shows a $10$dB improvement in SNR compared to that for the HSD estimator, and $15$ dB improvement in SNR compared to that for the Wiener Filter estimator.    

\revised{From the results in Fig.\,\ref{FigDParts21}, we can conclude that since the channel components in V2V channels (sparse and diffuse components) have different levels of energy, proximity operators such as SCAD and MCP with a multi-threshold nature of their proximal operators (see Fig.\,\ref{FigProximal} and Eq.s \eqref{multiSCAD} and \eqref{multiMCP}) are (more) suited to the channel structure.}
 
\begin{figure}
\centering
\includegraphics[width=2.8in]{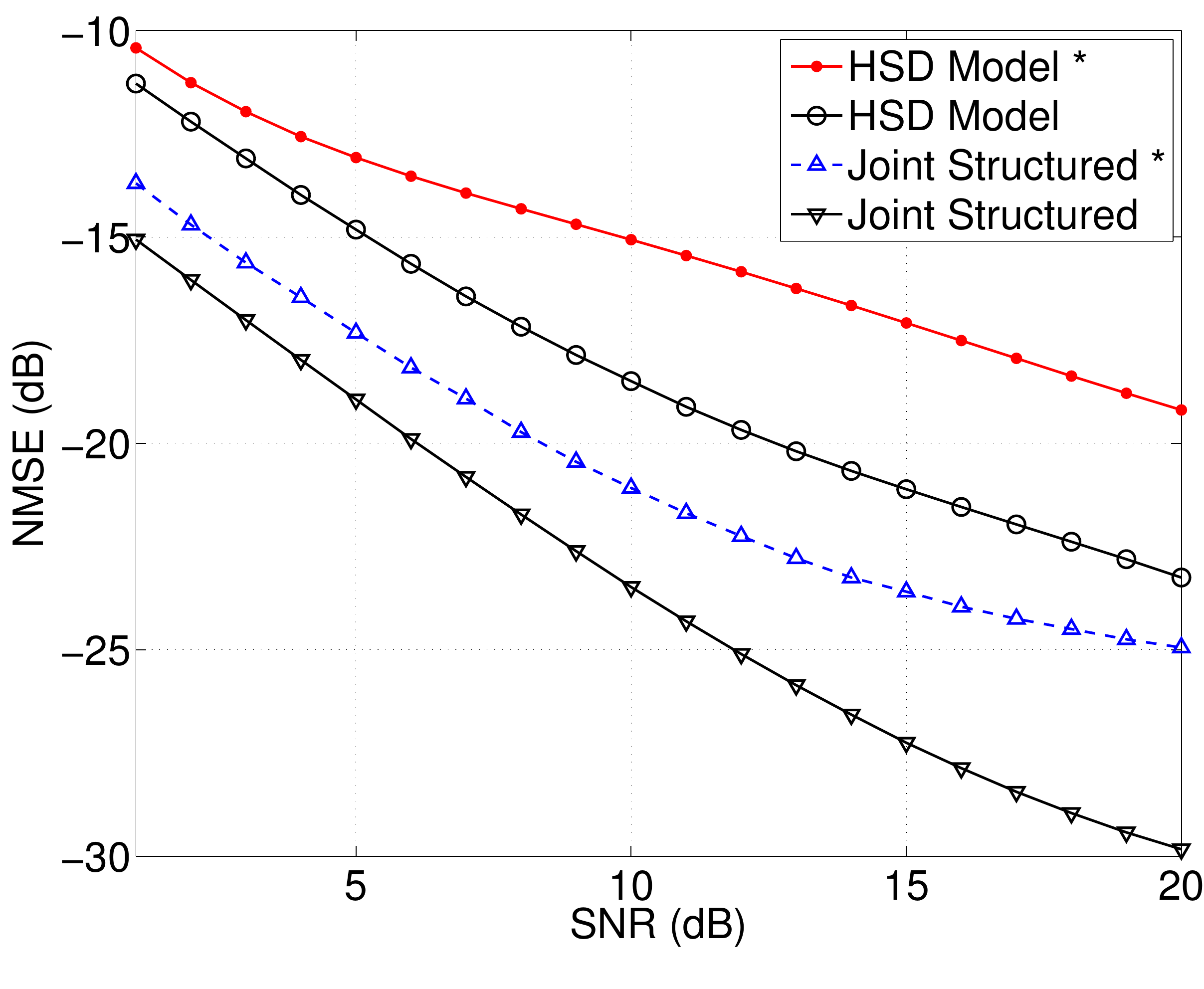}
\caption{NMSE of the channel estimators. The $(*)$ in the legend means that the leakage effect is not compensated, \emph{i.e.}, ${\bf G = I}$ is assumed in the channel estimation algorithm.} 
\label{FigDParts13}
\end{figure}

In Fig.~\ref{FigDParts13}, we consider the effect of leakage compensation (Section \ref{se:LekageComp}) on the performance of sparse estimator of V2V channels, such as the HSD estimator and our proposed joint sparsity estimator using SCAD regularizer function for group sparsity.  From the results in Fig.~\ref{FigDParts13}, we observe that the uncompensated leakage effect reduces performance severely, more than $7$ dB, particularly at higher SNR, due to the channel mismatch introduced by the channel leakage. 

\revised{
Next, we assess the performance of our proposed V2V channel estimation algorithm for different values of $N_{SD}$, $N_{DI}$, and $N_{MD}$ in the channel. In Fig.\,\ref{Sparsityorderofchan}.(a), the performance of our algorithm for  different choices of $N_{SD}+N_{MD}$ is depicted. Note that the static (SD) and mobile discrete (MD) components have similar effects on the channel sparsity pattern. Therefore, we consider the value of $N_{SD}+N_{MD}$ for our simulations. In Fig. \ref{Sparsityorderofchan}.(b), the performance of our algorithm for different choices of $N_{DI}$ is depicted
For this numerical simulation, we consider $K= 1024$, $M=512$, and $N_r = 2048$, and SNR=$10$ dB.
\begin{figure}[!t]
    \centering
    \subfloat[NMSE for varying $N_{SD}+N_{MD}$ and fixed $N_{DI}=500$]{{\includegraphics[width=0.3\textwidth]{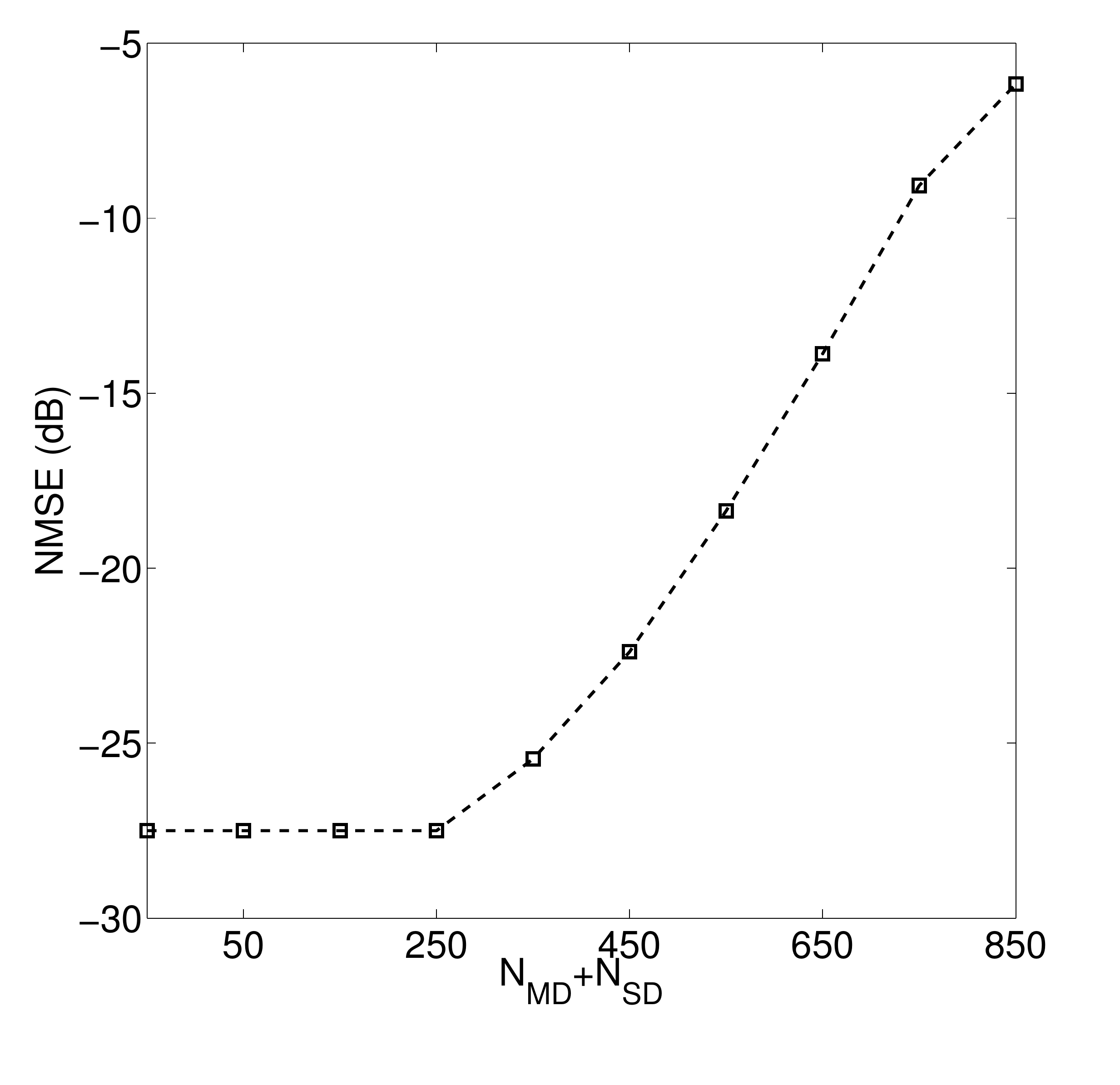} }}%
    \,
    \subfloat[NMSE for varying  $N_{DI}$  and fixed $N_{SD}+N_{MD}=50$]{{\includegraphics[width=0.3\textwidth]{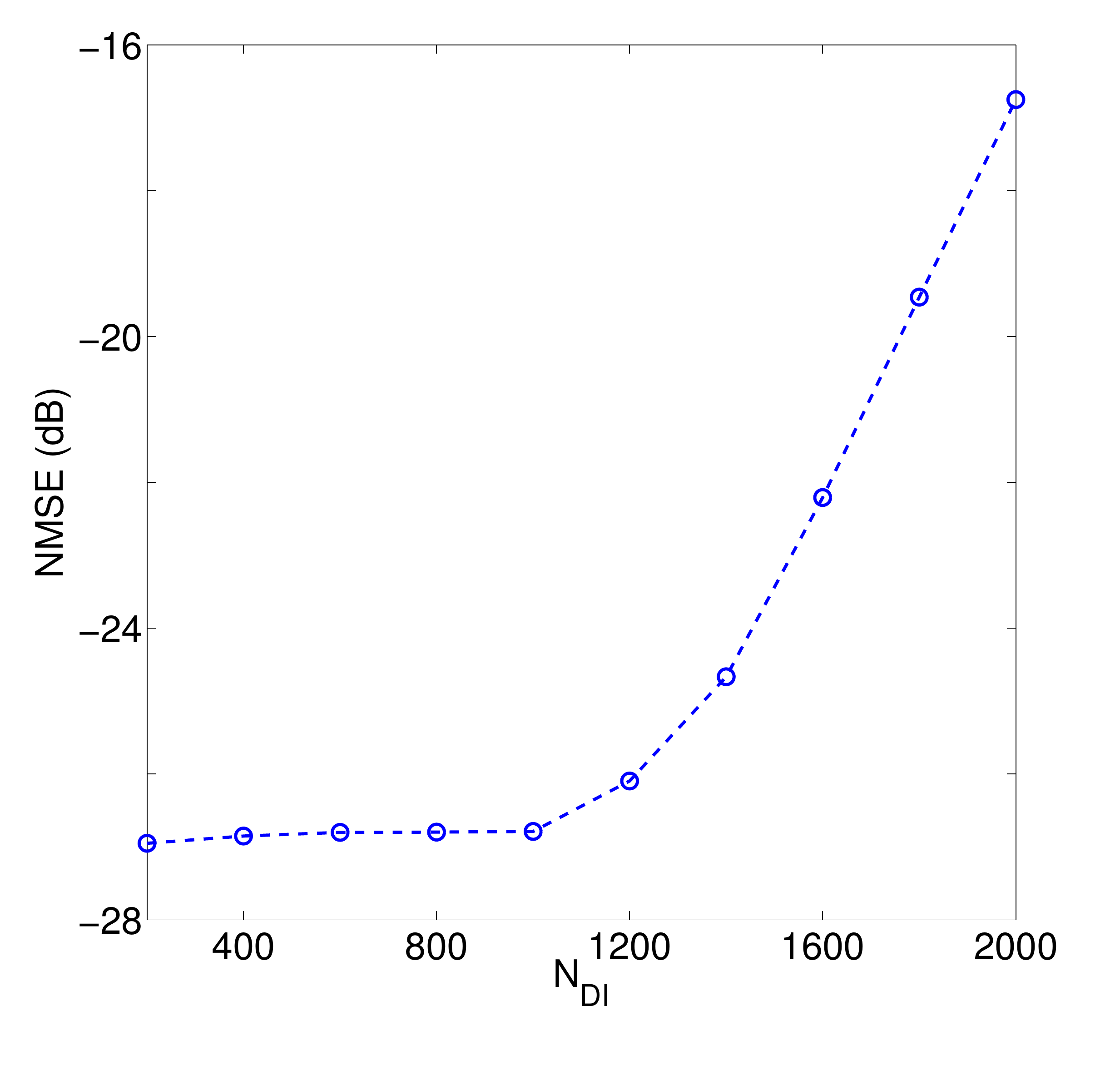} }}%
    \caption{Performance of algorithm for different values of  $N_{DI}$ and $N_{SD}+N_{MD}$.}%
\label{Sparsityorderofchan}
\end{figure}
It can be seen that by decreasing the number of channel components, the performance of our method improves.
}

\revised{Furthermore, we consider the performance of our proposed method under different values of $\Delta \tau$, and $\Delta \nu$. Our numerical results in Figures \ref{borderwidth}.(a) and \ref{borderwidth}.(b) show that by increasing both $\Delta\tau$ and $\Delta\nu$ the performance of our algorithm degrades. This is due to the fact that increasing $\Delta\tau$ and $\Delta\nu$ reduces the sparsity of the channel, which results in less noise reduction. Here $N_r = 512$, $M=256$, $N_{SD}+N_{MD}=20$ and $N_{DI}=200$ are considered.
\begin{figure}[!t]
    \centering
    \subfloat[Different choice of $\Delta\tau$]{{\includegraphics[width=0.3\textwidth]{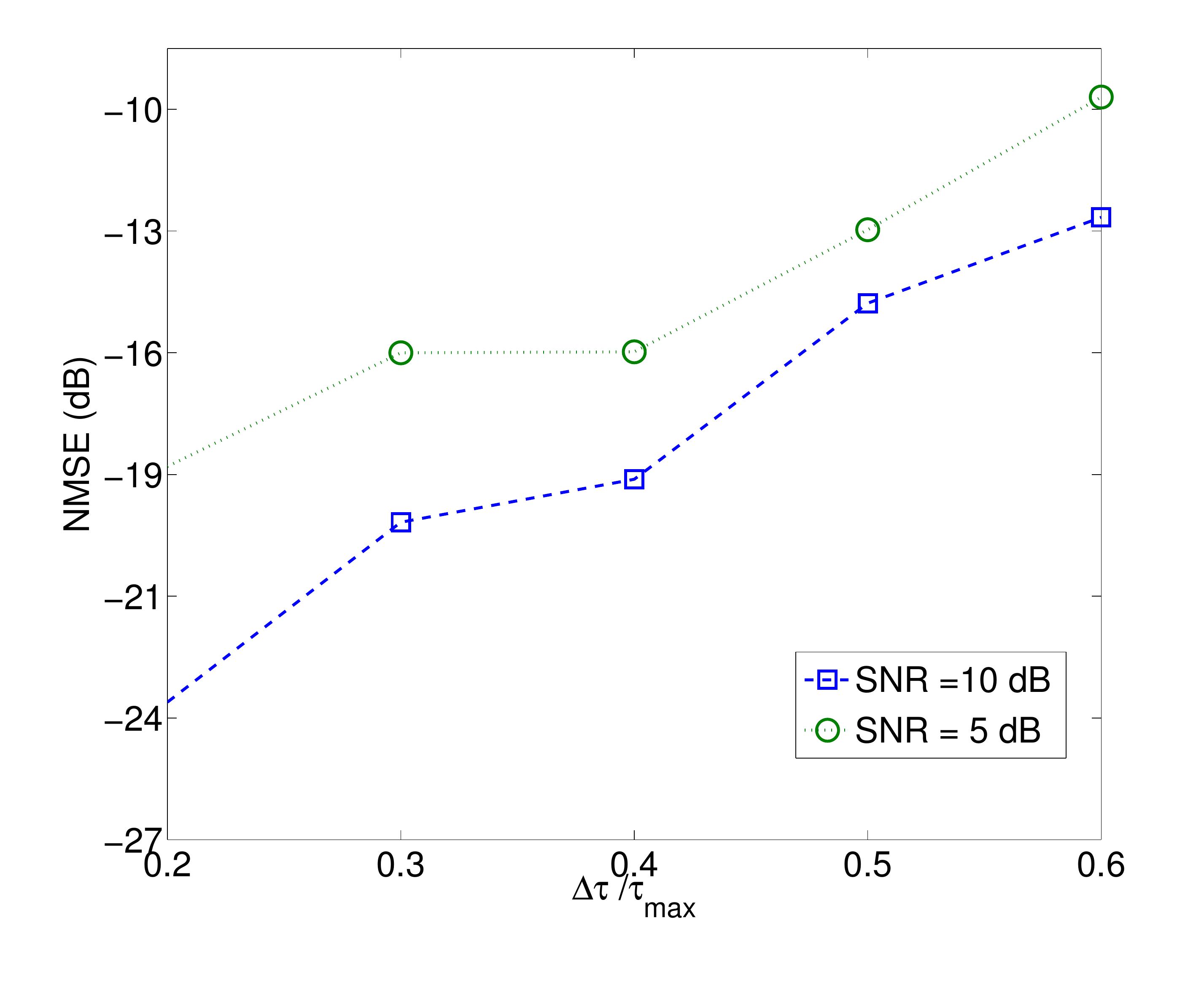} }}%
    \qquad
    \subfloat[Different choice of $\Delta\nu$]{{\includegraphics[width=0.3\textwidth]{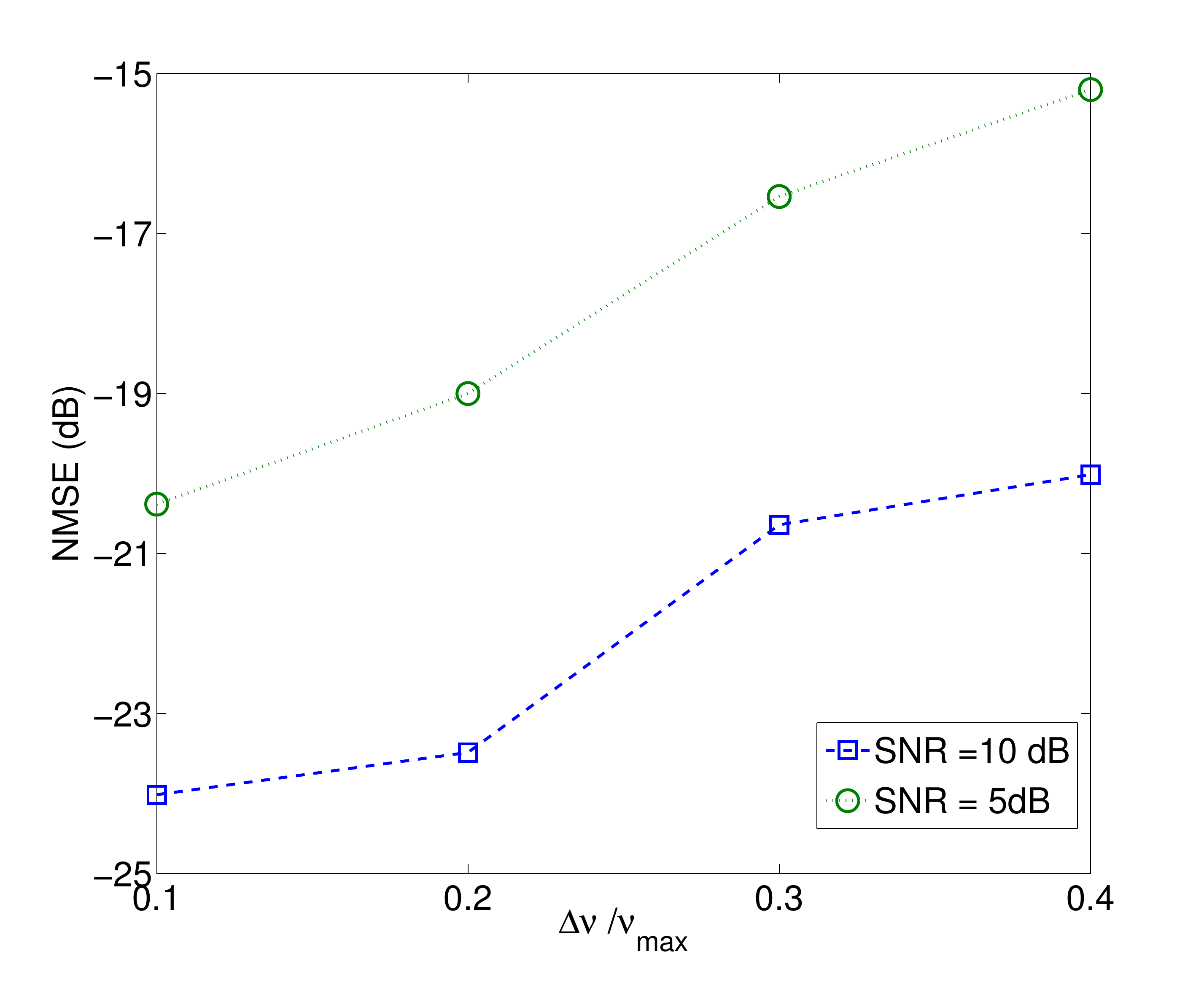} }}%
    \vspace{-0.1in}
    \caption{Performance of algorithm for different values of  $\Delta\tau$ and $\Delta\nu$.}%
\label{borderwidth}
\end{figure}
}

\revised{
The value of parameter $K$ is lower bounded by ${N_r}$, namely the total number of measurements, in the sense that $2K+1\ge N_r$, as given in Eq. \eqref{DDDRep}. To decrease the computational complexity, we consider $K=\lceil(N_r-1)/2\rceil$, which provides a good performance results as seen in Fig. \ref{Keffect}. 

However, we can consider $K \ge \lceil(N_r-1)/2\rceil$ by zero padding the DFT transform in Eq. \eqref{DDDRep}.
The larger the $K$, the better the leakage is compensated, but the higher the signal dimension of ${\bf x}$. The first will improve the signal recovery, but the latter degrades the signal reconstruction as more unknown variables are introduced. To understand the effect of increasing $K$, we consider $K= \gamma\lceil(N_r-1)/2\rceil$, where $1\le \gamma$. 
Here for computational efficiency, we consider  $\gamma = 2^n$, where $n$ is a non integer. Furthermore, $N_r = 512$, $M=256$, $N_{SD}+N_{MD}=20$ and $N_{DI}=200$ are considered.   Results in Fig.\ref{Keffect} show that by increasing $\gamma$, there is an improvement in the performance of the algorithm for $n\le4$.  But after $n\ge5$ (increasing the signal dimension) increasing $K$, increases the NMSE.

\begin{figure}[!t]
\centering
    \includegraphics[width=0.5\textwidth]{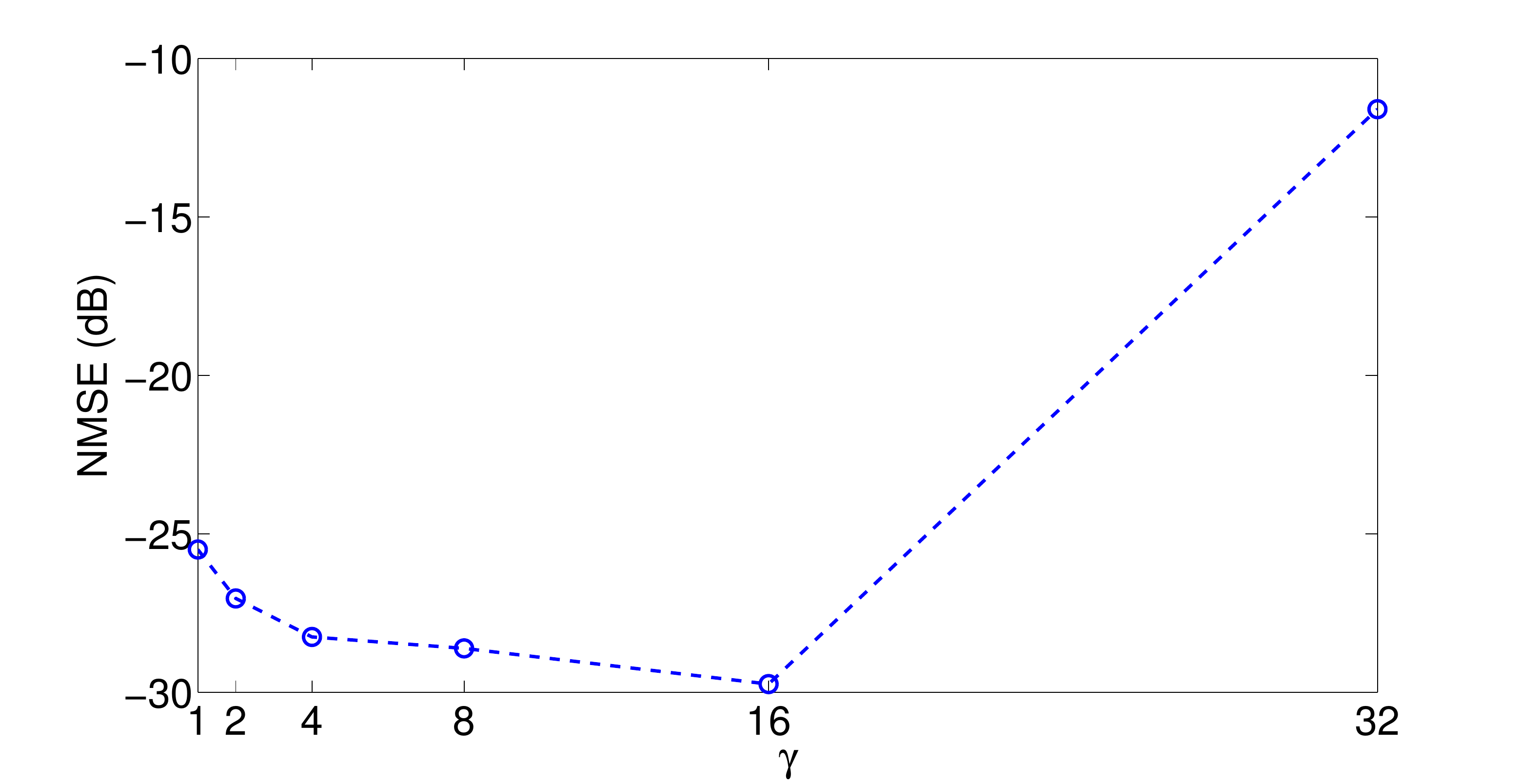}
    \caption{Performance of the proposed algorithm under different values of $K= \gamma\lceil(N_r-1)/2\rceil$.}
\label{Keffect}
\end{figure}}

\revised{In Figures \ref{Otherparam}.(a) and \ref{Otherparam}.(b), we have considered the effect of different choices of $M$ and $N_r$ on the NMSE of the proposed V2V channel estimation algorithm. Note that the values of $N_r$ and $M$ are dictated by the sampling time, training signal length and channel delay spread. Therefore, we can not increase the value of these parameters arbitrarily. Figures \ref{Otherparam}.(a) and \ref{Otherparam}.(b) indicate that increasing the number of measurements, $N_r$, and $M$, reduces the NMSE, as expected. 

\begin{figure}[!t]
    \centering
    \subfloat[Different choice of $N_r$ and $M=256$]{{\includegraphics[width=0.3\textwidth]{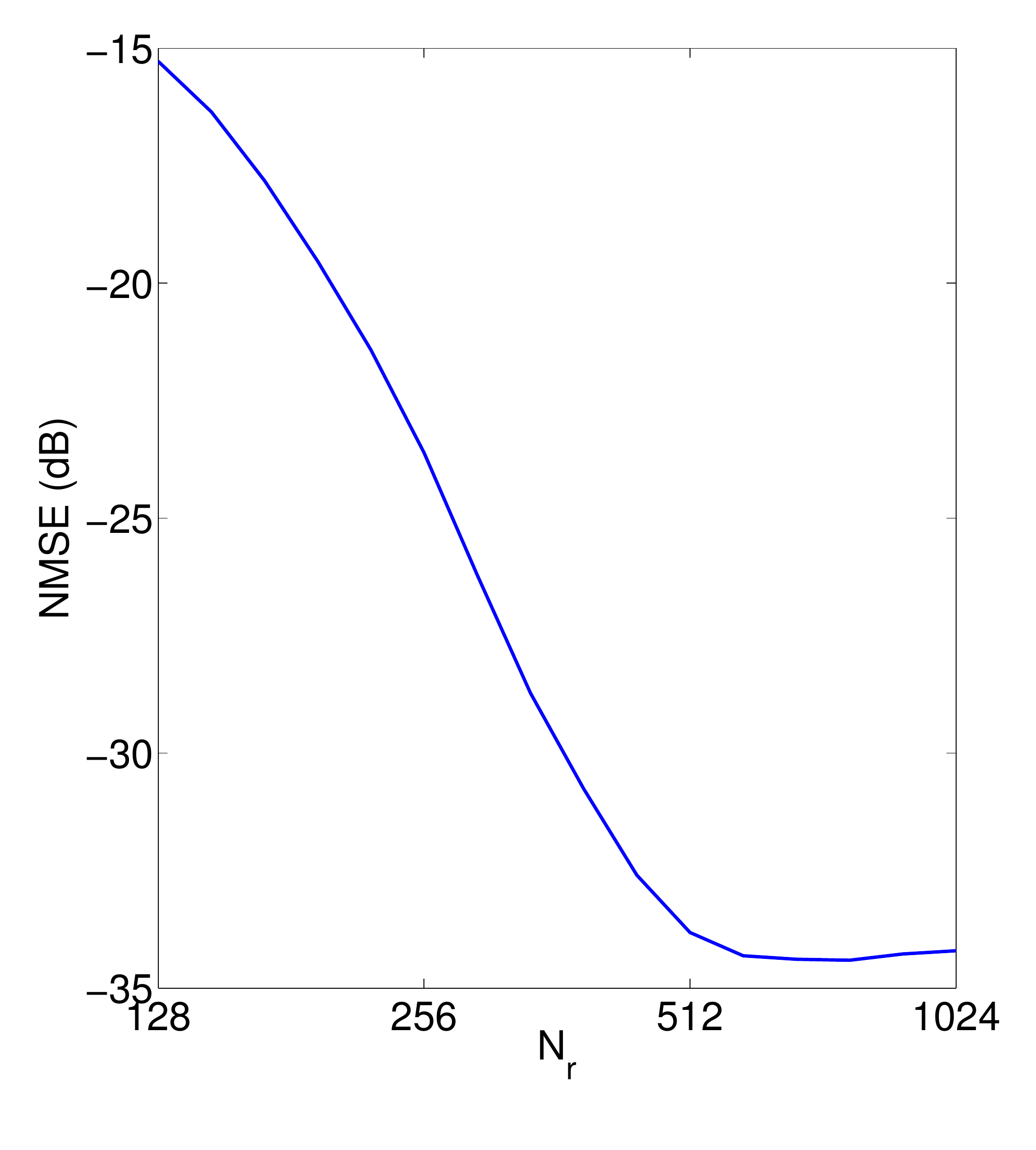} }}%
    \qquad
    \subfloat[Different choice of $M$ and $N_r = 256$]{{\includegraphics[width=0.3\textwidth]{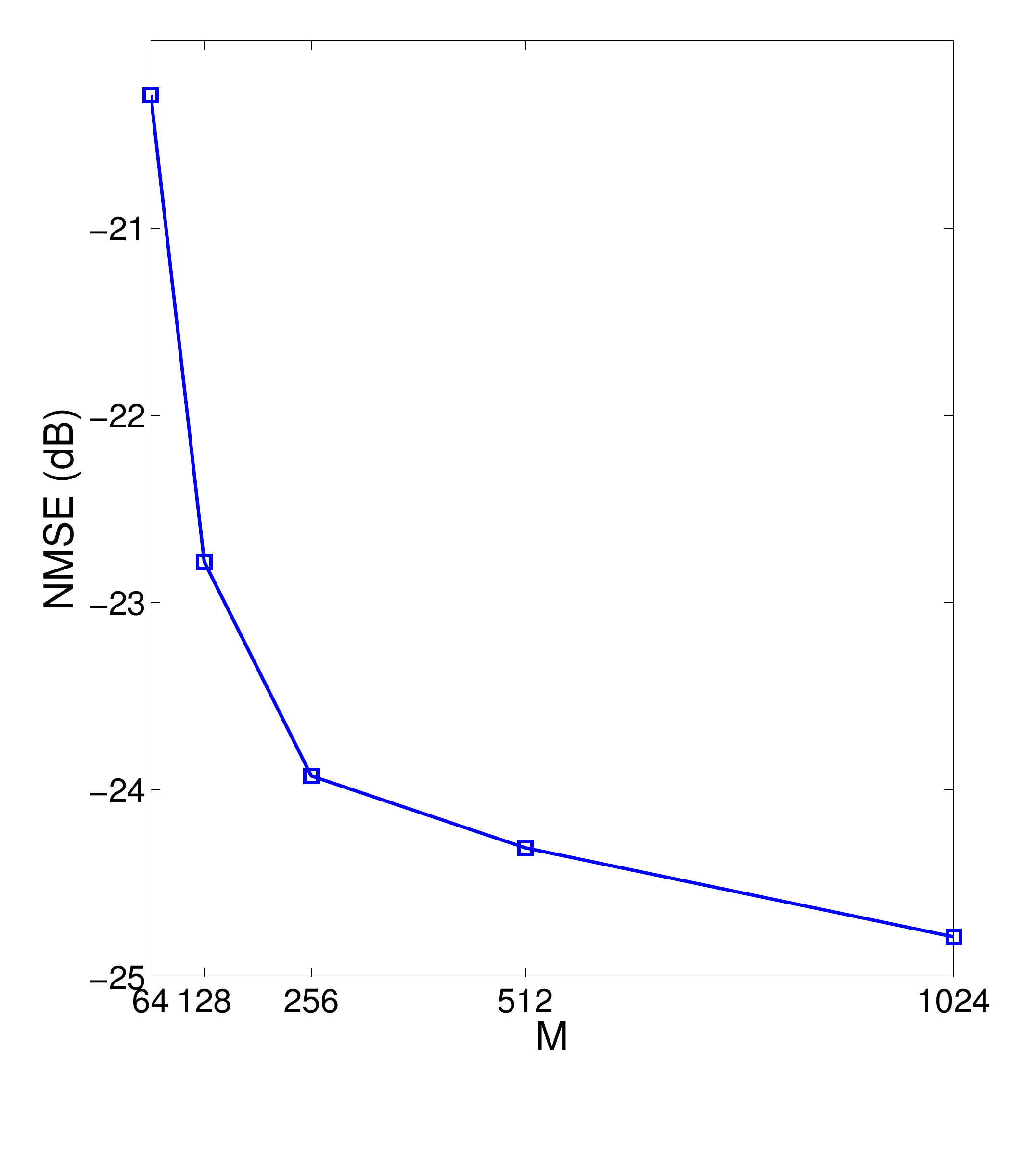}}}%
    \caption{Performance of the proposed algorithm under different values of $N_r$ and $M$.}%
    \label{Otherparam}%
\end{figure}}

\begin{figure}[!t]
\centering
    \includegraphics[width=0.4\textwidth]{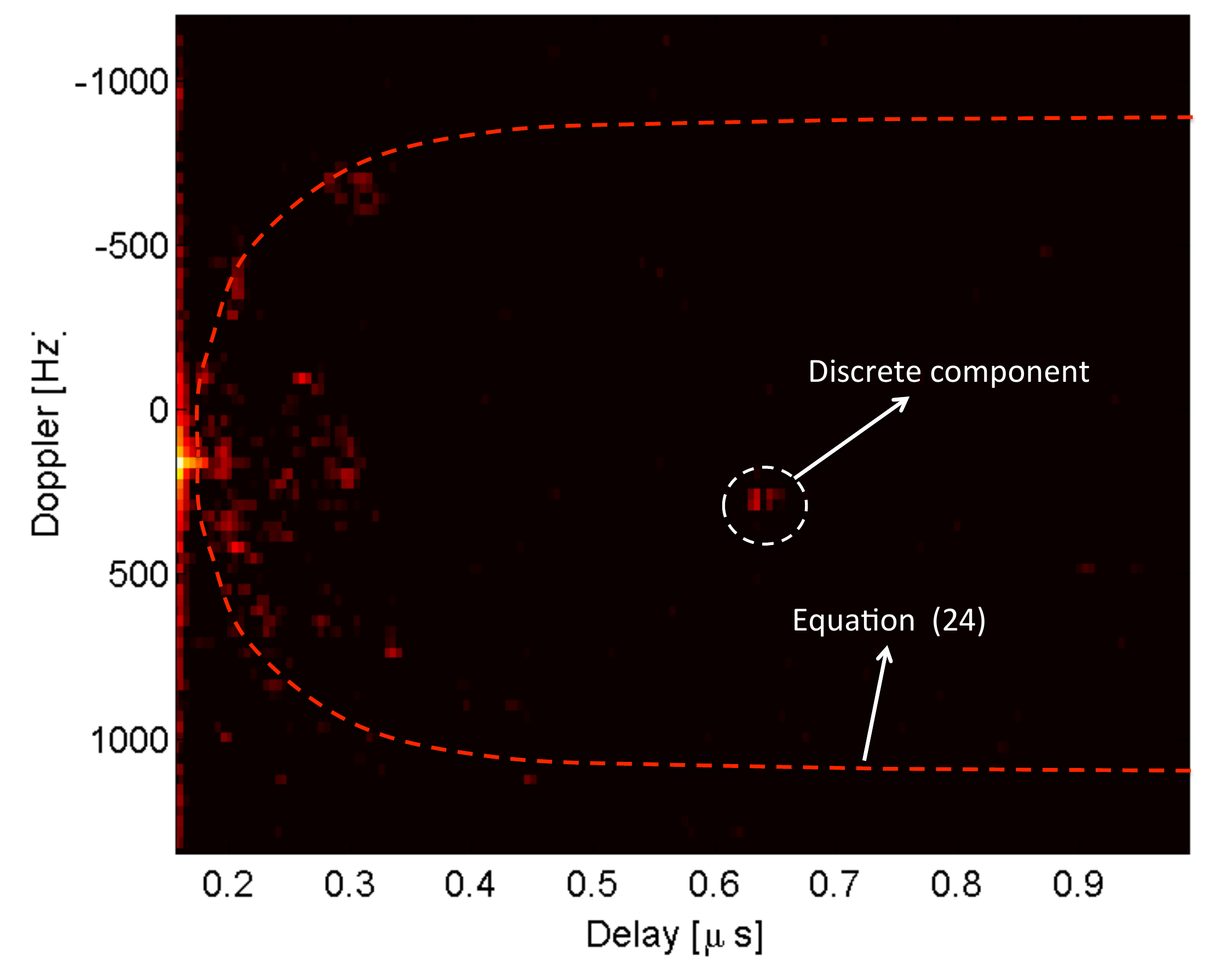}
    \caption{The channel delay-Doppler scattering function for a real channel measurement data \cite{Abbasthesis}.}
\label{PracticalFig0}
\end{figure}

\begin{figure}[!t]
    \centering
    \subfloat[$\Delta\tau$ = 0.3 $\mu$s and $\Delta\nu$ = 500 Hz]{{\includegraphics[width=0.4\textwidth]{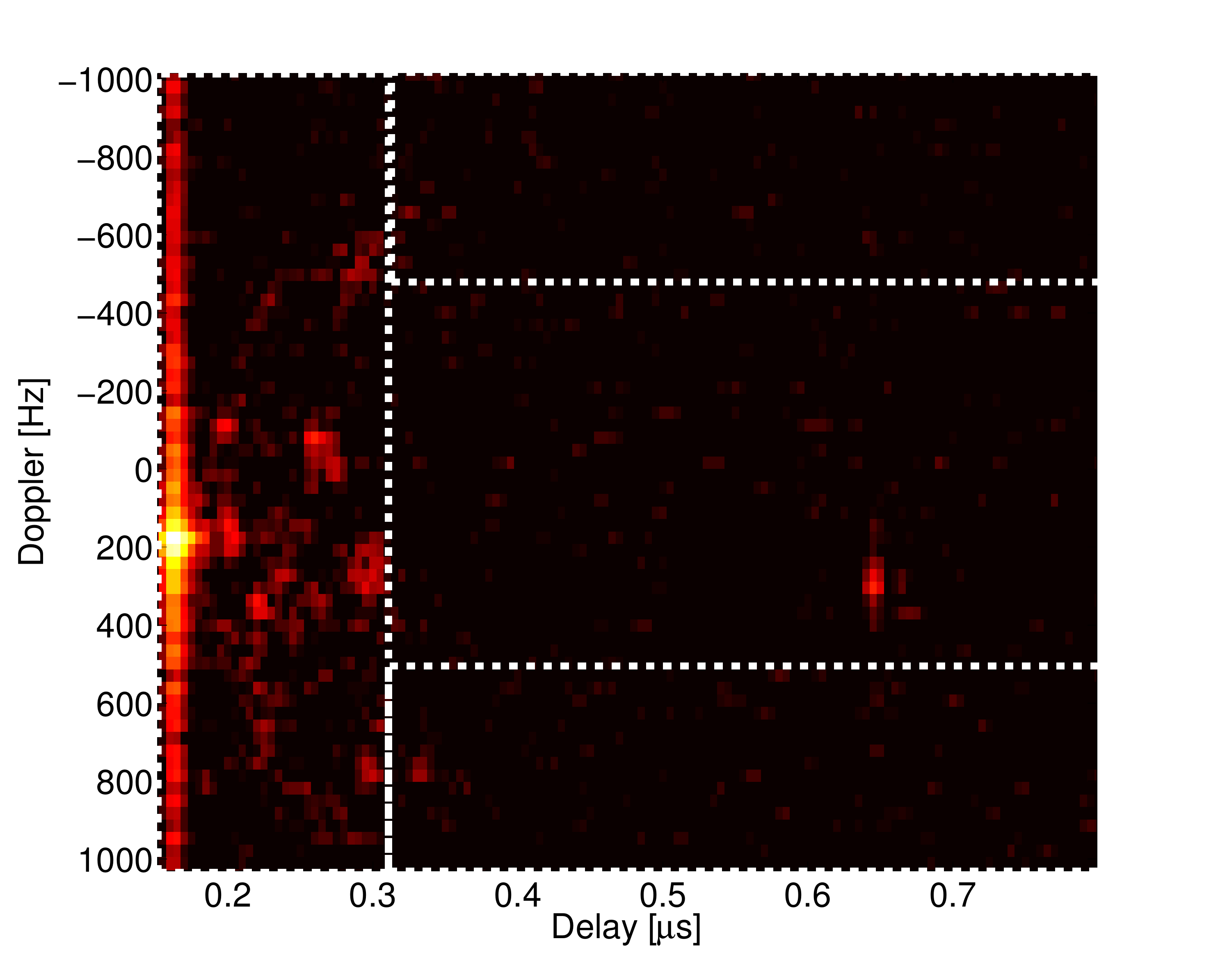} }}%
    \qquad
    \subfloat[$\Delta\tau$ = 0.33 $\mu$s and $\Delta\nu$ = 200 Hz]{{\includegraphics[width=0.4\textwidth]{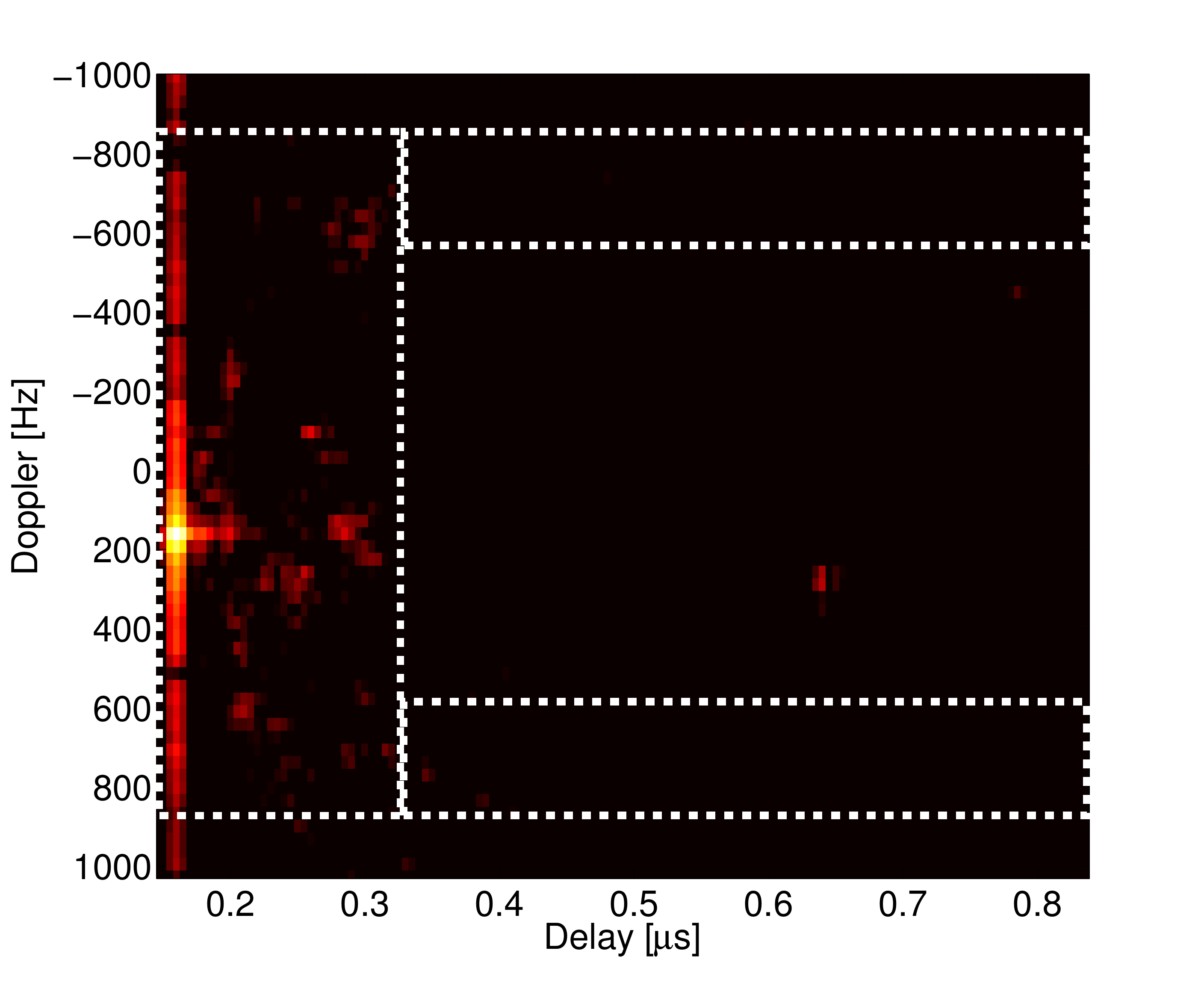}}}%
    \caption{Delay-Doppler spreading function. Diffuse components are confined to a U-shaped area.}%
    \label{PracticalFig34}%
\end{figure}

\begin{figure}[!t]
\centering
    \includegraphics[width=0.4\textwidth]{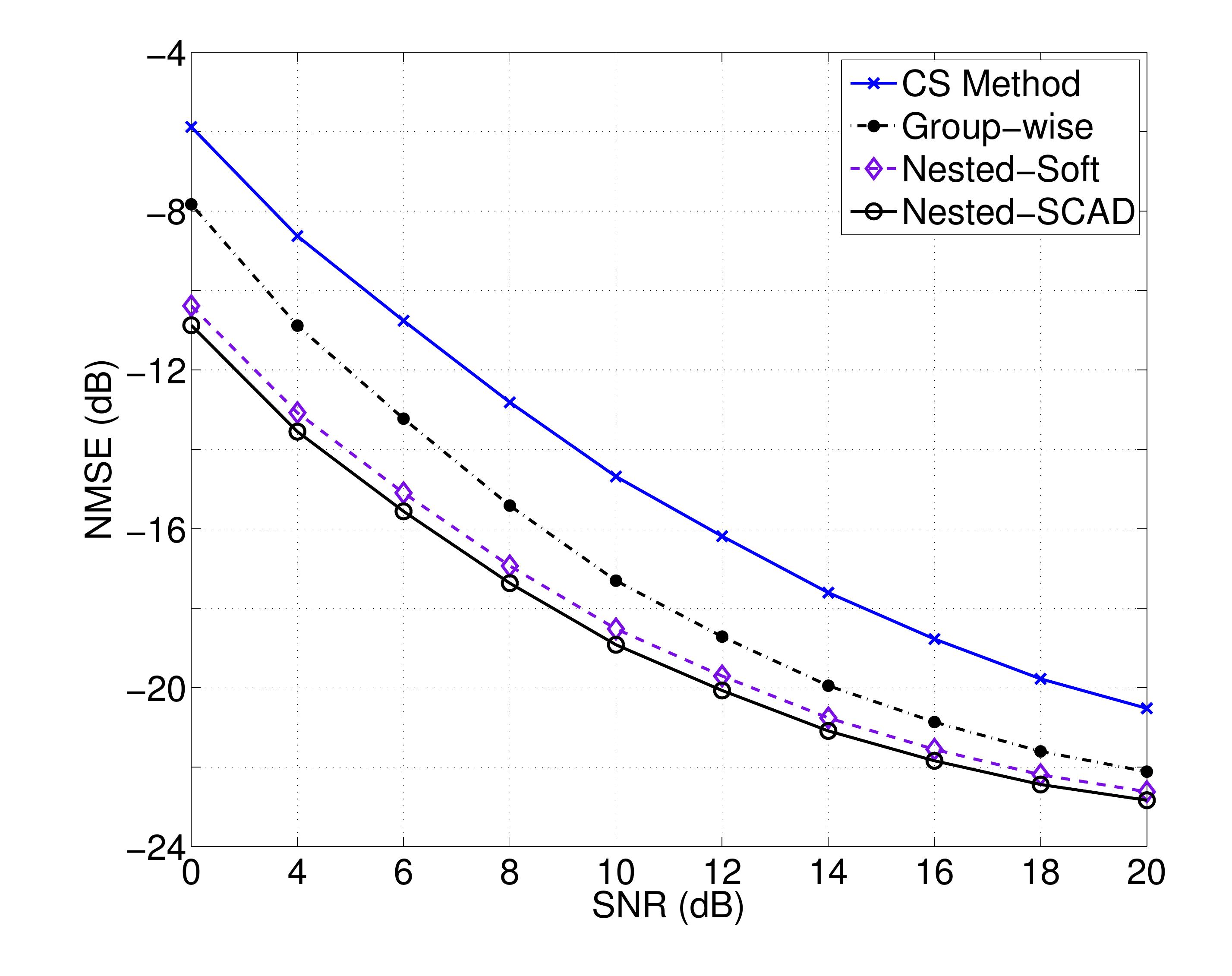}
    \caption{Comparison of NMSE v.s SNR for  CS Method \cite{bajwa2010compressed, taubock2010compressive}, and proposed method \emph{i.e.}, Group-wise ($\lambda_e=0$), Nested-soft \cite{Sprechmann2011CHiLasso, chartrand2013nonconvex},  and Nested-SCAD regularizers.}
\label{PracticalFig1}
\end{figure}

\begin{figure}[!t]
\centering
    \includegraphics[width=0.4\textwidth]{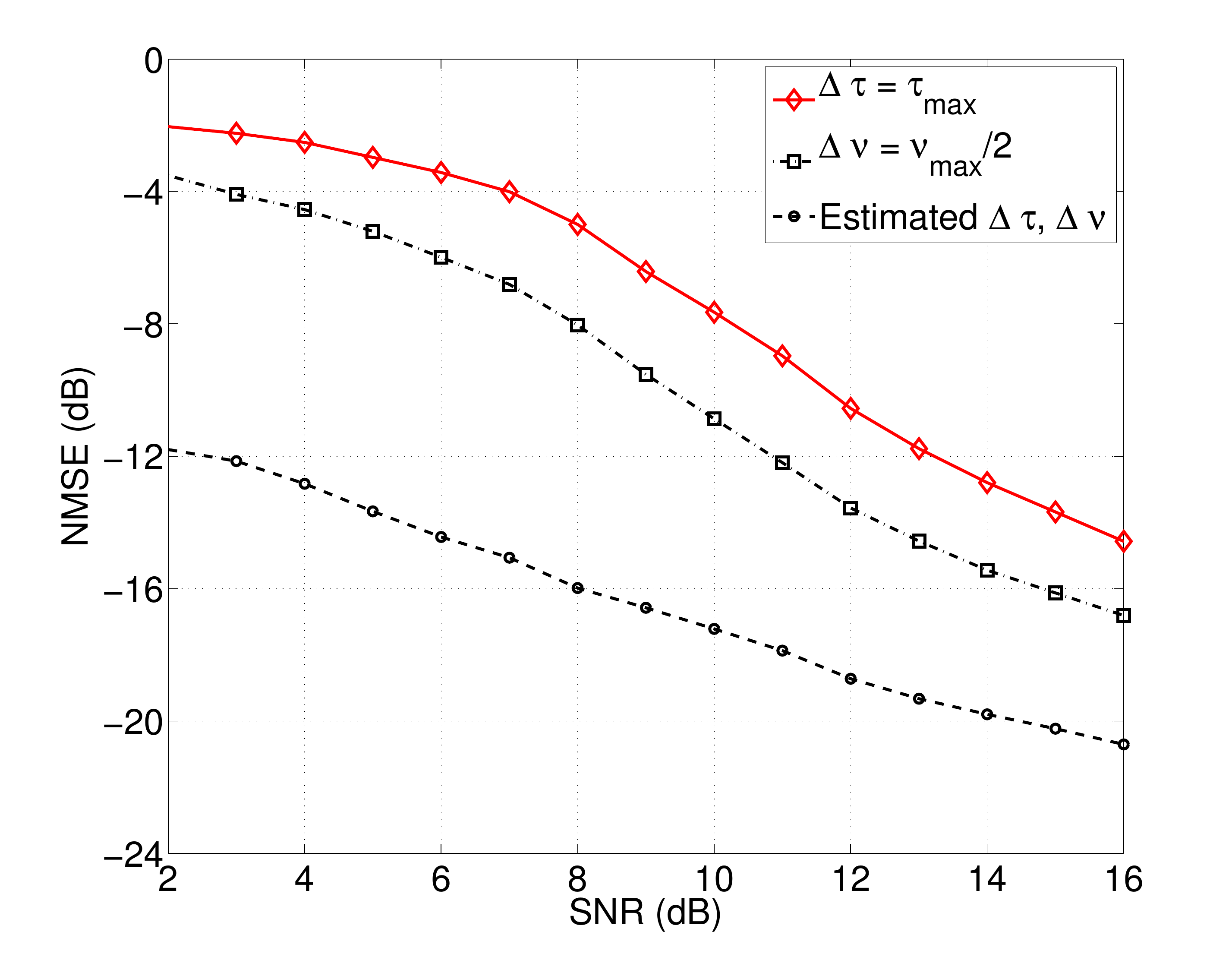}
    \caption{Performance of the proposed algorithm for different values of $\Delta\tau$
and $\Delta\nu$.}
\label{PracticalFig2}
\end{figure}

\revised{\section{Real Channel Measurements}\label{ExperimentSection}
In this section, we use experimental channel data, recorded by  WCAE\footnote{Wireless Communication in Automotive Environment} measurements  in a highway environment,  to model the V2V channel  and also assess the performance of our proposed channel estimation algorithm.  The channel measurement data was collected using the RUSK-Lund channel sounder in Lund and Malm\"o city in Sweden. The complex $4\times4$ multi-input-mutli-output (MIMO) channel transfer functions were recorded at  a $5.6$ GHz carrier frequency over a bandwidth of $240$ MHz in several different propagation environments with two standard Volvo V$70$ cars used as the TX and RX cars during the measurements \cite{Abbasthesis}. 

\subsection{Measurement setup}
V2V channel measurements were performed with the RUSK-LUND sounder using a multi-carrier signal with carrier frequency $5.6$ GHz to sound the channel and records the time-variant complex channel transfer function $H(t, f )$.  The measurement bandwidth was $240$ MHz, and a test signal length of $3.2 \mu$s was used.  The time-varying channel was sampled every $0.307$ ms, corresponding to a sampling frequency of $3255$ Hz during a time window of roughly $31.2$ ms. The sampling frequency implies a maximum resolvable Doppler shift of $1.5$ kHz, which corresponds to a relative speed of about $350$ km/h at $5.6$ GHz.

By Inverse Discrete Fourier Transforming (IDFT) the recorded frequency responses $H(t,f)$, with a Hanning window to suppress side lobes, the complex channel impulse responses $h(t, \tau)$ are obtained. Finally, by taking Discrete Fourier Transform (DFT) of $h(t, \tau)$ with respect to $t$,  the channel scattering function in the delay-Doppler domain,  $H[k,m]$, is computed. In this experiment, $K = 116$ and $M = 256$ is considered. 
Three recorded channel scattering functions, $H[k,m]$, are plotted in Figures \ref{PracticalFig0} and \ref{PracticalFig34}. 
Fig. \ref{PracticalFig0} presents a V2V channel in the delay-Doppler domain. A discrete component is visible at approximately $0.65\,\mu$s propagation delay. Also plotted in the figure is the Doppler shift vs. distance as produced by Equation \eqref{DopplerEq} and \eqref{delayEq}, \emph{i.e.}, for scatterers located on a line parallel to (and a distance $5$ m away from) the TX/RX direction of motion.  We notice that V2V channel scattering (in the delay-Doppler domain)  in Figures \ref{PracticalFig0} and  \ref{PracticalFig34} (a) and (b)  are highly structured as predicted in Section \ref{Ushape}. As seen in these figures,  the diffuse components are confined in a U-shaped area that was also predicted by our analysis in Section \ref{Ushape}. Furthermore,  we can observe the sparse structure of the discrete components in all the figures.   

In Fig. \ref{PracticalFig1}, we compare the performance of our proposed nested estimators with the CS method  \cite{bajwa2010compressed, taubock2010compressive} to estimate the channel given in Fig. \ref{PracticalFig34} (a).  The training pilot samples are generated as discussed in Section \ref{SimulationResult}. To vary the SNR, we add additive white Gaussian noise to the signal at the output of channel.  Note that bandwidth and observation time for the channel measurements are large enough to allow us to ignore the leakage effect. The results in Fig. \ref{PracticalFig1} confirm our numerical analysis in Section \ref{SimulationResult} and show that our proposed joint sparse estimation algorithm has a better performance compared to the (only) element-wise sparse estimator methods \cite{bajwa2010compressed, taubock2010compressive}. 

In Fig. \ref{PracticalFig2}, we investigate the effect of specifying regions on the channel estimation algorithm. We consider the channel given in Fig. \ref{PracticalFig34} (b) for this experiment. We  have considered three different region scenarios. In the first scenario, we determine the regions as computed by our heuristic method given in Appendix \ref{CoarseEstimation}. In the second scenario, we keep $\Delta \tau$ the same as the first scenario, but we set $\Delta \nu = \frac{\nu_{\rm max}}{2}$, which means that we neglect the structure of the diffuse components in regions $R_2$ and we assume that the diffuse components can occur in the entire delay-Doppler domain. Finally, in the third scenario, we  set $\Delta \tau= \tau_{\max}$, \emph{i.e.}, $R_1$ extends to cover the entire delay-Doppler domain.
Results in Fig. \ref{PracticalFig2} indicate that considering the structural information of the V2V channel in the delay-Doppler (three regions) significantly improves the performance of our joint sparse estimation algorithm. }
    
\section{Conclusions}
We provide a comprehensive analysis of V2V channels in the delay-Doppler domain using the well-known geometry-based stochastic channel modeling. Our characterization reveals that the V2V channel model has three key regions, and these regions exhibit different sparse/hybrid structures which can be exploited to improve channel estimation. Using this structure, we have proposed a joint element- and group-wise sparse approximation method using general regularization functions. We prove that for the needed optimization, the optimal solution results in a nested estimation of the channel vector based on the group and element wise penalty functions. Our proposed method exploits proximity operators and the alternating direction method of multipliers, resulting in a modest complexity approach with excellent performance. We characterized the leakage effect on the sparsity of the channel and robustified the channel estimator by explicitly compensating for pulse shape leakage at the receiver using the leakage matrix.
Simulation results reveal that exploiting the joint sparsity structure with non-convex regularizers yields a $5$ dB improvement in SNR compared to our previous state of the art HSD estimator in low SNR. Furthermore, using  experimental data of V2V channel from a WCAE measurement campaign, we showed  that our estimator yields $4$ dB to $6$ dB improvement in SNR compared to the compressed sensing method in \cite{bajwa2010compressed, taubock2010compressive}. 

\appendices

\section{Sparsity Inducing Regularizers}\label{RegularizerSectionVer}
 In this section, we show that the regularizer functions summarized in Section \ref{SUbsection111} satisfy Assumptions I.  We note that showing the assumptions for the convex regularizers is straightforward and thus omitted; we focus on the non-convex regularizers which will be used to induce group sparsity.\newline
{\bf SCAD regularizer \cite{fan2001variable}}: This penalty takes the form
\begin{align}
f_g(x;\lambda)=\begin{cases}\lambda |x|, & \text{for}\,\, |x|\le \lambda \\-\frac{x^2-2{\mu_S}\lambda |x|+\lambda^2}{2({\mu_S}-1)} & \text{for}\,\, \lambda<|x|\le {\mu_S}\lambda \\\frac{({\mu_S}+1)\lambda^2, }{2} &\text{for}\,\,|x|>{\mu_S}\lambda\end{cases}
\end{align}
where ${\mu_S}> 2$ is a fixed parameter. This penalty function is non-decreasing,  $f_g(0;\lambda)=f_g(x;0)=0$ and it is clear that $f_g(\alpha x; \alpha \lambda) =\alpha^2f_g(x;\lambda)$ for $\forall \alpha>0$. The derivative of the SCAD penalty function for $x\neq 0$ is given by
\begin{align}f'_g(x;\lambda)= {\rm sign}(x)\left(\frac{( {\mu_S}\lambda-|x|)_+}{{\mu_S}-1}\mathbb{I}\{|x|> \lambda\}+\lambda\mathbb{I}\{|x|\le \lambda\}\vphantom{\frac{( {\mu_S}\lambda-|x|)_+}{{\mu_S}-1}}\right) 
\end{align}
where $\mathbb{I}(.)$ is the indicator function, and any point in interval $[-\lambda, +\lambda]$ is a valid subgradient at $x=0$, so condition (iv) is satisfied. For $\mu = \frac{1}{{\mu_S}-1}$ the function $f_g(x;\lambda)+\mu x^2$ is also convex.
Thus, Assumption I holds for the SCAD penalty function with $\mu_S \ge 3$.\\
{\bf MCP regularizer \cite{zhang2010nearly}}: This penalty takes the form
\begin{align}
f_g(x;\lambda)={\rm sign}(x) \int_{0}^{|x|}\left(\lambda-\frac{z}{{\mu_M}} \right)_+\, dz
\end{align}
where ${\mu_M}> 0$ is a fixed parameter. This penalty function is non-decreasing for $x\ge 0$ and $f_g(0;\lambda)=f_g(x;0)=0$. Also,  $f_g(\alpha x; \alpha \lambda) =\alpha^2f_g(x;\lambda)$ for $\forall \alpha>0$. The derivative of the MCP penalty function for $x\neq 0$ is given by
\begin{align}
 f'_g(x;\lambda)=  {\rm sign}(x) \left(\lambda-\frac{|x|}{\mu_M} \right)_+
 \end{align}
and any point in $[-\lambda, +\lambda]$ is a valid subgradient at $x=0$. For $\mu = \frac{1}{\mu_M}$ the function $f_g(x;\lambda)+\mu x^2$ is also convex.
Thus, Assumption I holds for the MCP penalty function with $\mu_M\ge 2$.

\section{Proof of Theorem \ref{MainThm}}\label{The3}
We first prove two lemmas, needed for the proof of Theorem \ref{MainThm} and Corollary \ref{MainCoro}.

\begin{lem}\label{Lemmavec}
Consider that $g({\bf a};\lambda) = f\left(\frac{\|{\bf a}\|_2}{\rho};\lambda\right)$,
where $f(x;\lambda)$ is a non-decreasing function of $x$. Furthermore, $f(x;\lambda)$
is a homogeneous function, \emph{i.e.}, $f(\alpha x;\alpha\lambda) = \alpha^2f(x;\lambda)$.
Then,
\begin{enumerate}[i)]
  
\item $P_{\lambda, g}({\bf a}) = \gamma {\bf a}$, where $\gamma \in [0,1]$.
\item $\gamma  =\begin{cases}\frac{1}{\|{\bf a}\|_2}P_{\rho\lambda, \frac{f}{\rho^2}}(\|{\bf
a}\|_2) & \text{if}\,\, \| {\bf a} \|_2 > 0\\ 0 & \text{if}\,\, \| {\bf a} \|_2 = 0\end{cases}
$
\end{enumerate}
\end{lem}

{\bf Proof}: i). For every $\mat{z}$, we can write $\mat{z} = \mat{a}^{\perp}+\gamma
\mat{a}$, where $\mat{a}^{\perp} \perp \mat{a}$ and $\gamma\in \mathbb{R}.$
Therefore, based on the proximity operator definition we have,
\begin{align}
\label{HereIT}\prox{\lambda, g}(\mat{a}) &= \amin{\mat{z}}\LCr{\frac{1}{2}\Lp{\mat{a}-\mat{z}}{2}{2}+
g(\mat{z}; \lambda)}\\ \nonumber 
&= \amin{\mat{z} = \mat{a}^{\perp}+\gamma \mat{a}}\LCr{\frac{1}{2}\Lp{\mat{a}-\gamma\mat{a}-\mat{a}^{\perp}}{2}{2}+
g\LPr{\gamma \mat{a}+\mat{a}^{\perp};\lambda}}
\end{align}
Since $\Lp{\gamma\mat{a}+\mat{a}^{\perp}}{2}{}\ge \max\LCr{|\gamma|\|\mat{a}\|_2,
\|\mat{a}^{\perp}\|_2}$ and $f$ is a  non-decreasing function, we have $
g\LPr{\gamma \mat{a}+\mat{a}^{\perp};\lambda} \ge g\LPr{\gamma \mat{a};\lambda}.$ 
Therefore,  $\mat{a}^{\perp} = \mat{0}$ in the optimization problem \eqref{HereIT}
and we can rewrite it as,
\begin{align}\label{Proof2a}
\prox{\lambda, g}(\mat{a}) &=  {\mat{a}} \, \amin{\gamma}\LCr{\frac{1}{2}(\gamma-1)^2\|{\mat{a}}\|_{2}^{2}+
g\LPr{\gamma \mat{a};\lambda}}.
\end{align}
The two terms in the objective function in \eqref{Proof2a} are increasing
when $\gamma$ increases from  
$\gamma=1$ or decreases from $\gamma=0$. Hence, the minimizer lies in the
interval $[0, 1]$. Therefore, we have  $\prox{\lambda,
g}(\mat{a}) =\gamma \mat{a}$, where $\gamma \in [0,1]$.\\
ii). Let $t = \gamma \|\mat{a}\|_2$, then the optimization problem in
\eqref{Proof2a} for  $\| {\bf a} \|_2 > 0$, can be written as, 
\begin{align} \nonumber
 \prox{\lambda, g}(\mat{a})  &= \frac{\mat{a}}{\|\mat{a}\|_2}\,\,\amin{t\in
[0,\|\mat{a}\|_2]}\LCr{\frac{1}{2}\LPr{\|{\mat{a}}\|_{2}-t}^{2}+ f\left(\frac{t}{\rho};\lambda \right)}\overset{(a)}{=}\frac{\mat{a}}{\|\mat{a}\|_2}\,\, \amin{t\in [0,\|\mat{a}\|_2]}\LCr{\frac{1}{2}\LPr{\|{\mat{a}}\|_{2}-t}^{2}+
\frac{1}{\rho^2}f(t;\rho \lambda)}
 \\ 
&\overset{(b)}{=} \frac{\mat{a}}{\|\mat{a}\|_2}  \,\,P_{\rho\lambda,
\frac{f}{\rho^2}}(\|{\bf x}\|_2).
\end{align}
Equality (a) is due to the homogeneity of function $f$, \emph{i.e.},
$f\left(\frac{t}{\rho};\lambda\right) = f\left(\frac{t}{\rho};\frac{\rho \lambda}{\rho}\right)=
\frac{1}{\rho^2}f(t;\rho \lambda)$. Equality (b) is due to the definition
of the proximal operator.  Thus, 
 $\gamma =  \frac{1}{\|{\bf a}\|_2}P_{\rho\lambda, \frac{f}{\rho^2}}(\|{\bf
a}\|_2),$ and the proof of the Lemma is completed.

\begin{lem}\label{Homegity}
If the function $f(x;\lambda)$ is homogenous \emph{i.e.}, $f(\alpha x;\alpha\lambda)=\alpha^2f(x;\lambda)$
for all $\alpha>0$, then $ \prox{\alpha\lambda, f}(\alpha{b}) = \alpha\prox{\lambda, f}({b}) $ 
for $\forall {b}\in \RR{}$ and $\lambda>0$. 
\end{lem}
{\bf Proof}: By definition of the proximity operator, we have $P_{{\alpha}\lambda, f}({\alpha}b) = \amin{{x}}\LCr{\frac{1}{2}\LPr{{{\alpha}b}-x}^{2}+ f(x;{\alpha}\lambda)}$. Consider $x={\alpha}z$, and using the homogenous properties of $f$, we have $$P_{{\alpha}\lambda, f}({\alpha}b) = {\alpha}\,\amin{{z}}\LCr{\frac{{\alpha}^2}{2}\LPr{{b}-z}^{2}+ {\alpha}^2f(z;\lambda)}={\alpha}\,\amin{{z}}\LCr{\frac{1}{2}\LPr{{b}-z}^{2}+f(z;\lambda)}={\alpha}P_{\lambda, f}(b).$$

{\bf Proof of Theorem \ref{MainThm}:} 
Since the regularizer functions $\phi_e$ and $\phi_g$ are separable, it is easy to show that the solution of optimization problem in Equation \eqref{KeyOpti} can be computed in parallel for all the groups as, 
$ {\hat{\bf{a}}}_i  =  \underset{{\bf a}_i \in \mathbb{R}^{N}}{\operatorname{argmin}} 
 \left\{\frac{1}{2}\left\| {\bf{b}}_i-{\bf{a}}_i\right\|_2^2+g({\bf a}_i;\lambda_{g})+E({\bf a}_i;\lambda_{e})\right\}$
for $i=1, \dots, N_g$, where $g({\bf a}_i;\lambda_{g})  =  f_g \left({\|{\bf a}_i \|_2} ; \lambda_{g} \right)$
and $ E({\bf a}_i;\lambda_{e}) =  \sum_{j }  f_e \left( {{\bf a}_i[j]};  \lambda_{e} \right)$.
For the sake of simplicity in notation of the proof for Theorem \ref{MainThm} and  Corollary \ref{MainCoro},  we drop the group  index and we consider
${\hat{\bf{a}}}  =  \underset{{\bf a}}{\operatorname{argmin}} 
 \left\{\frac{1}{2}\left\| {\bf{b}}-{\bf{a}}\right\|_2^2+g({\bf a};\lambda_{\rho g})+E({\bf a}_i;\lambda_{\rho e})\right\}$
where $g({\bf a};\lambda_{\rho g})=f_g\left(\frac{\|{\bf{a}}\|_2}{\rho};\lambda_{\rho
g}\right)$ and $E({\bf a};\lambda_{\rho e}) = \sum_{j}f_e\left(\frac{a[j]}{\rho};\lambda_{\rho
e}\right)$. Here $\lambda_{\rho g} = \frac{\lambda_{g}}{\rho}$ and $\lambda_{\rho e} = \frac{\lambda_{g}}{\rho}$. Note that for $\rho =1$, we have the claim in Theorem \ref{MainThm}. 

Assume $\mat{v} = \prox{\lambda_{\rho e}, E}(\mat{b})$ and $\mat{u} = \prox{\lambda_{\rho
g}, g}(\mat{v})$. 

Based on above definitions, to prove the claim of Theorem \ref{MainThm},
we need to show that $\mat{a} = \mat{u}$ is the minimizer of $J(\mat{a})=\frac{1}{2}\left\|{\bf b} - {\bf a}\right\|_2^2+g({\bf a};\lambda_{\rho g})+E({\bf a};\lambda_{\rho e})$.
To prove this claim, we consider two cases: 
I:  $\mat{u} \neq \mat{0}$, and II: $\mat{u}= \mat{0}$.\\

Case (I): $\mat{u} \neq \mat{0}$.  Since $\mat{u} = \prox{\lambda_{\rho g},
g}(\mat{v})$ and $g({\bf a};\lambda_{\rho g})=f_g\left(\frac{\|{\bf{a}}\|_2}{\rho};\lambda_{\rho g}\right)$, and $f_g$ is a homogenous non-decreasing function, by Lemma \ref{Lemmavec}
we have  $\mat{u} = \gamma \mat{v}$, for some $\gamma \in (0,1]$. Furthermore,
 $\mat{u}$ should satisfy the first order optimality condition for the objective
function in $\mat{u} =\amin{\mat{a}}\LCr{\frac{1}{2}\Lp{\mat{v}-\mat{a}}{2}{2}+ g(\mat{a};
\lambda_{\rho g})},$ namely
\begin{align}\label{E1}
\mat{0} \in \mat{u}-\mat{v}+ \subg{g}(\mat{u};\lambda_{\rho g}).
\end{align}
Using the definition of the proximity operator and Remark \ref{rem22}, we have $[\prox{\lambda_{\rho
e}, E}(\mat{b})]_i = \prox{\lambda_{\rho e}, f_e}\left(\frac{{b}[i]}{\rho}\right)$.
Since $f_e$ is a homogeneous function, using Lemma \ref{Homegity}, we have
$\prox{\gamma\lambda_{\rho e}, f_e}\left(\gamma\frac{{b}[i]}{\rho}\right)=
\gamma\prox{\lambda_{\rho e}, f_e}\left(\frac{{b}[i]}{\rho}\right)
$ or equivalently, 
\begin{align}\nonumber
\prox{\gamma\lambda_{\rho e}, E}(\gamma \mat{b}) &= \amin{\mat{a}}\LCr{\frac{1}{2}\Lp{\gamma\mat{b}-\mat{a}}{2}{2}+
E(\mat{a}; \gamma\lambda_{\rho e})}
=\gamma\,\amin{\mat{z}}\LCr{\frac{1}{2}\Lp{\mat{b}-\mat{z}}{2}{2}+ E(\mat{z};
\lambda_{\rho e})} \\&\label{thisonehere0} =\gamma\prox{\lambda_{\rho e}, E}( \mat{b}) =\gamma \mat{v}=\mat{u}
\end{align}
and by the first order optimality condition (of $\mat{u}$) for the objective
function in \eqref{thisonehere0}, we have 
$\mat{0} \in \mat{u}-\mat{\gamma b}+\subg{E(\mat{u}; \gamma\lambda_{\rho
e})}.$
Since $\gamma \neq 0$, above Equation can be rewritten as 
\begin{align}\label{inamig}
\mat{0} \in \mat{v}-\mat{b}+\frac{1}{\gamma}\subg{E({\bf u}; \gamma\lambda_{\rho
e})}.
\end{align} 
Since $E({\bf u};\lambda_{\rho e}) = \sum_{j}f_e\left(\frac{u[j]}{\rho};\lambda_{\rho
e}\right)$, applying scale invariant property of function $f_e$, \emph{i.e.}, $f_e\left(\frac{u[j]}{\rho};\gamma\lambda_{\rho
e}\right)  = \gamma f_e\left(\frac{u[j]}{\rho};\lambda_{\rho e}\right)$,
we have $\frac{1}{\gamma}\subg{E({\bf u}; \gamma\lambda_{\rho e})}=\subg{E({\bf
u}; \lambda_{\rho e})}$. Therefore, we can rewrite \eqref{inamig} as
\begin{align}\label{E2}
\mat{0} \in \mat{v}-\mat{b}+\subg{E({\bf u};\lambda_{\rho e})}.
\end{align} 
Summing Equations \eqref{E1} and \eqref{E2}, we have 
$\mat{0} \in \mat{u}-\mat{b}+\subg{g(\mat{u}; \lambda_{\rho g} )}+ \subg{E(\mat{u};\lambda_{\rho
e})},$ which is the first order optimality of $\mat{u}$ for the objective function
$J(\mat{a})$.
Case (II): $\mat{u} = 0$. Here, we need to show the first-order optimality conditions for $\mat{u} = 0$ for the objective function $J(\mat{a})$, \emph{i.e.}, $\mat{0} \in \LCr{\mat{u}-\mat{b}+\subg{g(\mat{u}; \lambda_{\rho
g} )}+ \subg{E(\mat{u};\lambda_{\rho e})}}|_{\mat{u}=\mat{0}}$ $= \subg{g(\mat{0};
\lambda_{\rho g} )}+ \subg{E(\mat{0};\lambda_{\rho e})}-\mat{b}.$
This is equivalent to showing the existence of a $\mats{\chi}_1 \in [-1,+1]^{N}$,
 \big(equivalent to the term $ \subg{E(\mat{0};\lambda_{\rho e})}$\big),
where  $\mats{\chi}_2$ with $\|\mats{\chi}_2\|_2\le
\frac{1}{\rho}$ \big(equivalent to the term $\subg{g(\mat{0}; \lambda_{\rho g} )}$\big)
such that $\mat{b} = \lambda_{\rho e} \mats{\chi}_1 + \lambda_{\rho e}
\mats{\chi}_2$, due to property (iii) in Assumption I. By definition of the proximity operator we have 
\begin{align} \label{FO2132}
\mat{u} = \prox{\lambda_{\rho g},
g}(\mat{v}) &= \amin{\mat{z}}\LCr{\frac{1}{2}\Lp{\mat{v}-\mat{z}}{2}{2}+ g(\mat{z};
\lambda_{\rho g})} = \amin{\mat{z}}\LCr{\frac{1}{2}\Lp{\mat{v}-\mat{z}}{2}{2}+ f_g\left(\frac{\|{\bf{z}}\|_2}{\rho};\lambda_{\rho g}\right)}.
\end{align}
Using the first optimality condition of $\mat{u} = \mat{0}$ for the objective function in \eqref{FO2132}, we have 
\begin{align}\label{E1001}
\mat{0} \in -\mat{v}+ \subg{f_g}(0;\lambda_{\rho g})\frac{\subg{(\|\mat{0}\|)}}{\rho}.
\end{align}
Since  $\subg{(\|\mat{0}\|_2)} = \left\{\mat{x} \in \mathbb{R}^{N}, \|\mat{x}\|_2 \le 1\right\}$ and $|z| \le \lambda_{\rho g}$ for all $z \in \subg{f_g}(0;\lambda_{\rho g})$ (using property (iii) in Assumption I),  using Equation \eqref{E1001} we have 
$\|\mat{v}\|_2 \le \frac{\lambda_{\rho g}}{\rho}$. 
Furthermore, since $\mat{v} = \prox{\lambda_{\rho e}, E}(\mat{b})$,  the first-order optimality
condition implies that $\mat{0} \in  \subg{E(\mat{v}; \lambda_{\rho e})}+\mat{v}-\mat{b}$.
Thus for $\mats{\chi}_1 \in \subg{E(\mat{v}; \lambda_{\rho e})}$ and $\mats{\chi}_2 = \frac{\mat{v}}{\lambda_{\rho g}}$, we have $\mat{b} = \lambda_{\rho e} \mats{\chi}_1 + \lambda_{\rho e}
\mats{\chi}_2$ and proof is completed.

\section{Region and Group Specification}\label{CoarseEstimation}
\begin{figure}[t]
\centering
\includegraphics[width=2.5in]{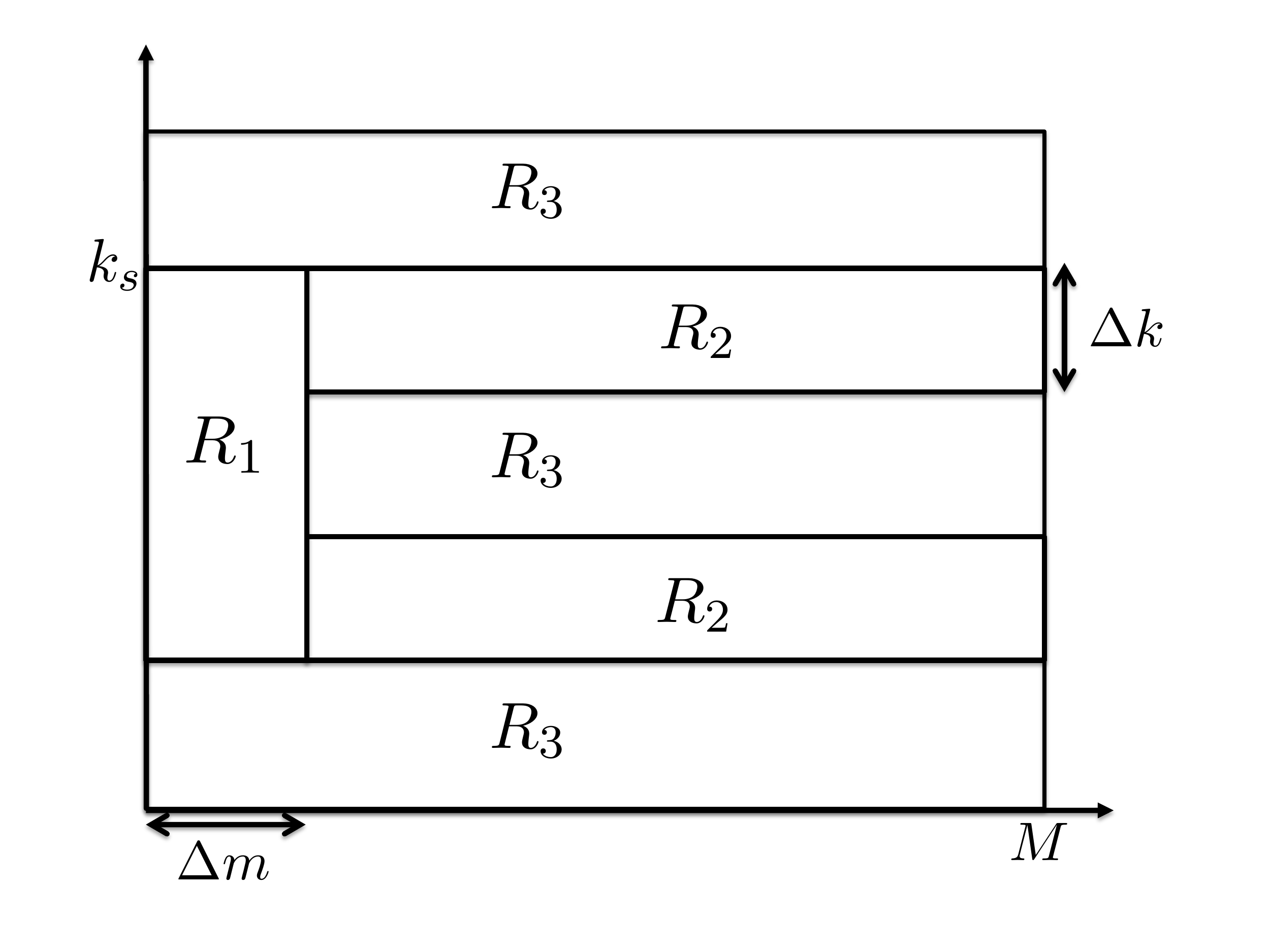} 
\caption{Discrete representation of the regions $R_1$, $R_2$, and $R_3$.}
\label{fig:regions}
\end{figure}

Here, we propose a heuristic method to find the regions $R_1$,  $R_2$, and $R_3$,  depicted in Fig.~\ref{fig:regions}, and introduced in Section ~\ref{Ushape}.  To describe the regions, we need to compute the value of $k_S$, $\Delta k$, and $\Delta m$, \emph{i.e.}, the discrete Doppler and delay parameters that corresponds to $\nu_{S}$, $\Delta\nu$, and $\Delta\tau$, respectively. To estimate the  $\Delta m$ and $\Delta k$, we use a regularized least-squares estimate of
$\mathbf{x}$ given by  ${\bf{x}}_{LS}={\bf A}_0{\bf{y}}={\bf A}_0 (\mathbf{A} {\bf{x}} + \mathbf{z}) \approx {\bf x}+ {\bf e}$ where ${\bf e}={\bf{A}}_0{\bf z}$ and ${\bf A}_0=({\bf A}^H{\bf A}+\rho^2 {\bf I})^{-1}{\bf A}^H$ and $\rho$ is a small real value. Based on relationship in the Eq.~\eqref{Vectorizing}, we can write
${\bf{x}}_{LS} = \vstack \left\{{\bf H}_{LS}\right\}$, where ${\bf H}_{LS}$ is an estimate of the discrete delay-Doppler spreading function.
Let us define the function $\mathcal{E}_d \left(m\right) = \frac{1}{m}\sum_{j=1}^{m}\sum_{i=-K}^{+K}|H_{LS}[i,j]|^2$ for $1\le m \le M$. This function represent the energy profile of channel components in delay direction. Then,
\begin{align}
\Delta m = \min_m  \left\{m | \mathcal{E}_d \left(m\right) \le T_1 \right\}, 
\end{align}
where $T_1= \alpha_d \mathcal{E}_{d}\left(0\right)$. Here $0<\alpha_d < 1$ is a tuning parameter. For a highway environment \cite{karedal2009geometry}, based on our numerical analysis,  $\alpha_d \approx 0.4$ is a reasonable choice.

After computing $\Delta m$, to compute the values of  $\Delta k$ and $k_s$ as labeled in Fig.\,8, 
due to symmetry of channel diffuse components around zero Doppler value, we define functions $\mathcal{E}_{\nu}\left(k\right)=\sum_{i=\Delta m}^{M}|H_{LS}[k,i]|^2+|H_{LS}[-k,i]|^2$ for $0\le k \le K$. This function represents the energy profile of channel components in Doppler direction. Let us define $k_0 = \max_k{\mathcal{E}_{\nu}\left(k\right)}$. Then, we can estimate $\Delta k$ and $k_s$ using following equations,
\begin{align}
k_s-\Delta k &= \max_k \left\{k |  \mathcal{E}_{\nu}\left(k\right)<T_2,  0\le k < k_0 \right\},\\
k_s &= \min_k \left\{k |  \mathcal{E}_{\nu}\left(k\right)<T_2,  k_0< k \le K \right\},
\end{align}
where $T_{\nu} = \alpha_{\nu}\mathcal{E}_{\nu}\left(k_0\right)$ with $0<\alpha_{\nu}< 1$.
For highway environment, based on our numerical analysis,  $\alpha_{\nu} \approx 0.6$ is a good choice.
\section{Proximal ADMM Iteration Development}\label{gabay1976dualAPEN}
For the optimization problem given in \eqref{OptimizationMain}, we form the augmented Lagrangian
\begin{align}
L_{\rho}\left({\bf{x}},{\bf{w}}, {\boldsymbol{\theta}}\right) &= \frac{1}{2}\left\| {\bf{y}}-{\bf{A}}{\bf{x}}\right\|_2^2+ \phi_g(|{\bf{w}}|;\lambda_g)+ \phi_e(|{\bf{w}}|;\lambda_e)+\left<{\boldsymbol{\theta}}, {\bf{x}}-{\bf{w}}\right> +\frac{\rho^2}{2}\|{\bf{x}}-{\bf{w}}\|_2^2,
\end{align}
where ${\boldsymbol{\theta}}$ is the dual variable, $\rho\neq0$ is the augmented Lagrangian parameter, and $\left<{\bf a}, {\bf b}\right> = \text{Re}({\bf b}^{H}{\bf a})$. Thus, ADMM consists of the iterations:
\begin{itemize}
\item {\bf update}-{\bf x}:
${{\bf{x}}}^{n+1} = \underset{{\bf{x}}}{\operatorname{argmin}}\,\frac{1}{2}\left\| {\bf{y}}-{\bf{A}}{\bf{x}}\right\|_2^2 +\left<{\boldsymbol{\theta}}^{n} , {\bf{x}}-{\bf{w}}^{n}\right> +\frac{\rho^2}{2}\|{\bf{x}}-{\bf{w}}^{n} \|_2^2. $

\item {\bf update}-{\bf w}:
${\bf{w}}^{n+1} = \underset{{\bf{w}}}{\operatorname{argmin}}\,\phi_g(|{\bf{w}}|;\lambda_g)+ \phi_e(|{\bf{w}}|;\lambda_e) +\left<{\boldsymbol{\theta}}^{n} , {\bf{x}}^{n+1} -{\bf{w}}\right>+\frac{\rho^2}{2}\|{\bf{x}}^{n+1} -{\bf{w}}\|_2^2.
$
\item {\bf update-dual variable}: ${\boldsymbol{\theta}}^{n+1} =  {\boldsymbol{\theta}}^{n} + \rho^2 \left({\bf{x}}^{n+1}-{\bf{w}}^{n+1}\right).$
\end{itemize}
Deriving a closed form expressions for update-${\bf x}$ is straightforward,
${\bf{x}}^{n+1} = \rho^2{\bf{A}}_0 \left({\bf{w}}^{n} -{\boldsymbol \theta}_{\rho}^{n}\right) + {\bf{x}}_0,$
where ${\boldsymbol \theta}_{\rho}^{n} = \frac{{\boldsymbol \theta}^{n}}{\rho^2}$, ${\bf{A}}_0 = \left(\rho^2{\bf{I}}+{\bf{A}}^{H}{\bf{A}}\right)^{-1}$ and ${\bf{x}}_0 =  {\bf{A}}_0{\bf{A}}^H{\bf{y}}$. 
If we pull the linear terms into the quadratic ones in the objective function of update-${\bf w}$ and ignoring additive terms, independent of ${\bf w}$, then we can express this step as
\begin{align}\nonumber
{\bf{w}}^{n+1} =  \underset{{\bf w}}{\operatorname{argmin}} \left\{\frac{1}{2}\left\|{\bf{x}}^{n+1}+{\boldsymbol \theta}_{\rho}^{n} - {\bf w}\right\|_2^2 + \frac{1}{\rho^2}\left(\phi_g(|{\bf{w}}|;\lambda_g) 
+ \phi_e(|{\bf{w}}|;\lambda_e)\right) \vphantom{\frac{1}{2}\left\|{\bf{x}}^{n+1}+{\boldsymbol \theta}_{\rho}^{n} - {\bf w}\right\|_2^2} \right\}\\
= \underset{{\bf w}}{\operatorname{argmin}} \sum_{i=1}^{N_g}\left\{\frac{1}{2}\left\| {\bf{x}}^{n+1}_i+{\boldsymbol \theta}_{\rho i}^{n} - {\bf w}_i\right\|_2^2 
 +f_g\left(\frac{\|{\bf{w}}_i\|_2}{\rho};\frac{\lambda_g}{\rho}\right)
 +\sum_{j}f_e\left(\frac{|w_i[j]|}{\rho};\frac{\lambda_e}{\rho}\right) \vphantom{\frac{1}{2}\left\| {\bf{x}}^{n+1}_i+{\boldsymbol \theta}_{\rho i}^{n} - {\bf w}_i\right\|_2^2} \right\}
\end{align} 
where ${\bf x}_i$, ${\bf w}_i $,  and ${\boldsymbol \theta}_{\rho i}$ are computed using the partitions introduced for the channel vector in Section \ref{ITmeansHere}. 
Thus, we can perform the update-${\bf w}$ step in parallel for all groups,
\begin{align}
{\bf{w}}^{n+1}_i =  \underset{{\bf w}_i}{\operatorname{argmin}} \left\{\frac{1}{2}\left\| {\bf{x}}^{n+1}_i+{\boldsymbol \theta}_{\rho i}^{n} - {\bf w}_i\right\|_2^2 
 +f_g\left(\frac{\|{\bf{w}}_i\|_2}{\rho};\frac{\lambda_g}{\rho}\right)
 +\sum_{j}f_e\left(\frac{|w_i[j]|}{\rho};\frac{\lambda_e}{\rho}\right) \vphantom{\frac{1}{2}\left\| {\bf{x}}^{n+1}_i+{\boldsymbol \theta}_{\rho i}^{n} - {\bf w}_i\right\|_2^2}\right\}
\end{align}
Here, for simplicity in representation, we define $\lambda_{\rho g}=\frac{\lambda_g}{\rho}$ and $\lambda_{\rho e}=\frac{\lambda_e}{\rho}$. 
In addition, we define  $E(|{\bf w}_i|;\lambda_{\rho e}) = \sum_{j}f_e\left(\frac{|w_i[j]|}{\rho};\lambda_{\rho e}\right),$ $g(|{\bf w}_i|;\lambda_{\rho g}) = f_g\left(\frac{\|{\bf{w}}_i\|_2}{\rho};\lambda_{\rho g}\right).$
Thus, we have
\begin{align}
{\bf{w}}^{n+1}_i &=  \underset{{\bf w}_i}{\operatorname{argmin}} \left\{\frac{1}{2}\left\|{\bf{x}}_i^{n+1}+{\boldsymbol \theta}_{\rho i}^{n} - {\bf w}_i\right\|_2^2 
 +g(|{\bf w}_i|;\lambda_{\rho g}) +E(|{\bf w}_i|;\lambda_{\rho e})  \vphantom{\frac{1}{2}\left\|{\bf{x}}_i^{n+1}+{\boldsymbol \theta}_{\rho i}^{n} - {\bf w}_i\right\|_2^2} \right\}
\end{align} 

To guarantee convergence to the optimal solution in $(P_0)$, the overall objective function,  \emph{i.e.}, $\frac{1}{2}\left\| {\bf{y}}-{\bf{A}}{\bf{x}}\right\|_2^2+ \phi_g(|{\bf{x}}|;\lambda_g)+ \phi_e(|{\bf{x}}|;\lambda_e)$,  should be a convex function \cite{eckstein2012augmented}. Note that since the first term in the objective function, \emph{i.e.}, the quadratic penalty function, is convex, for any functions $f_e$ and $f_g$ that satisfies condition (iv) in Assumption I, the overall objective function is convex as well. Thus, ADMM yields convergence for all choices of the convex and non-convex functions given in Section \ref{SUbsection111}.

\section{Proof of Lemma \ref{LemComplexCOnv}}\label{ProofLemma}
 The function $\|\mathbf{c} - \mathbf{w}\|_2^2 
= \|\mathbf{c}\|_2^2 + \|\mathbf{w}\|_2^2-2\text{Re}\{\mathbf{c}^H\mathbf{w}\}
= \|\, |\mathbf{c}|\, \|_2^2 + \| \,|\mathbf{w}|\, \|_2^2-2\text{Re}\{\mathbf{c}^H\mathbf{w}\}
$ is minimized, with respect to phase of $\mat{w}$, when $\text{Re}\{\mathbf{c}^H\mathbf{w}\}$ is maximized. Now,
\begin{align}
  \text{Re}\{\mathbf{c}^H\mathbf{w}\} = \sum_{n=1}^N |\mathbf{c}[n]| |\mathbf{w}[n]|\cos({\rm Ang}(\mathbf{c}[n])-{\rm Ang}(\mathbf{w}[n]))
\le \sum_{n=1}^N |\mathbf{c}[n]| |\mathbf{w}[n]| = |\mathbf{c}|^T |\mathbf{w}|
\end{align}
with equality if and only if ${\rm Phase}(\mathbf{w}) = {\rm Phase}(\mathbf{c})$, which in turn implies that $\|\mathbf{c} - \mathbf{w}\|_2^2 = \|\,|\mathbf{c}| - |\mathbf{w}|\,\|_2^2$. Hence,
\begin{align}
\amin{|\mathbf{w}| \odot {\rm Phase}(\mathbf{w}) \in \mathbb{C}^N}
\|\mathbf{c} - \mathbf{w}\|_2^2 = \amin{|\mathbf{w}| \odot {\rm Phase}(\mathbf{c}) \in \mathbb{C}^N}
\|\,|\mathbf{c}| - |\mathbf{w}|\,\|_2^2= {\rm Phase}(\mathbf{c})\odot\amin{\mathbf{|w|}\in \mathbb{R}^N} \|
\,|\mathbf{c}| - |\mathbf{w}|\,\|_2^2  
\end{align}
and the lemma follows. 

\section*{Acknowledgments}
We would like to thank Fredrik Tufvesson and Taimoor Abbas of Lund University for providing channel measurement data. The experimental work was funded by the Swedish Governmental Agency for Innovation Systems - VINNOVA, through the project Wireless Communication in Automotive Environment.  

\bibliographystyle{IEEEtran}
\bibliography{MyRefV2V.bib}

\end{document}